\documentclass[a4paper,10pt]{article}
\usepackage[cp1251]{inputenc}
\usepackage[left=44mm,right=40mm,top=60mm,bottom=52mm]{geometry}

\usepackage{amsmath}
\usepackage{amssymb}

\newtheorem{theorem}{Theorem}

\newtheorem{corollary}{Corollary}
\newtheorem{definition}{Definition}
\newtheorem{lemma}{Lemma}
\newtheorem{proposition}{Proposition}
\newtheorem{remark}{Remark}
\usepackage{mathrsfs} 
\begin{document}
\section*{About some extended Erlang-Sevastyanov queueing system and its convergence rate.}
\begin{flushright}
\textbf{G. A. ZVERKINA}\\

{\it Moscow State University of \\Railway Engineering (MIIT),\\
V.A. Trapeznikov Institute of Control\\
Sciences of the Russian\\
Academy of Sciences\\
zverkina@gmail.com\\}
\end{flushright}
UDC 519.21\\
{\small \textbf{Key words :} distribution of the queueing system state, rate of convergence, regenerative processes,
coupling method.
\begin{center}
\textbf{Abstract}
\end{center}
The upper bound for the convergence rate of the distribution  of a queuing system state with infinitely many servers is obtained for the case when the intensities of the incoming  and service flows depend on the state of the system.}

\section{Introduction. }
\subsection{Motivation}

The queueing system with infinitely many servers is considered. 
Let $t_1$, $t_2$,\ldots, $t_n$ be the times of first, second,\ldots , $n$-th customers input (denote $t_0\stackrel{\text{\rm def}}{=\!=} 0$). 
The customer immediately begins to be serviced at the input time. 
The  the time of the service of $i$-th customer is a random variable (r.v.) -- $\xi_i$ with cumulative distribution function (c.d.f) $G_i(s)$. The periods between customer inputs are $\eta_i \stackrel{\text{\rm def}}{=\!=} t_i-t_{i-1}$ -- r.v. with c.d.f. $F_i(s)$. For simplicity, suppose that the distributions of r.v. $\xi_i$ and $\eta_i $ are absolutely continuous (with respect to the Lebesgue measure). Denote $f_i(s)\stackrel{\text{\rm
def}}{=\!=} F'_i(s)$, $g_i(s)\stackrel{\text{\rm def}}{=\!=}
G'_i(s)$.

The distributions $F_i(s)$ and $G_i(s)$ can be described by the {\it intensities}.
\begin{definition}\label{zv1} 
Let $\mathscr{T}$ be a some random time period (with c.d.f. $\Phi(x)$) starting at the time $t=0$, and  $\mathscr{T}$ is not ended let at the time $s>0$. 
Then $\mathsf{P} \{\mathscr{T} \in (s,s+\Delta]| \mathscr{T}>s\}=\displaystyle \frac{\Phi(s+\Delta)-\Phi(s)}{1-\Phi(s)}$.

If c.d.f. $\Phi(x)$ is absolutely continuous, then
$$
\displaystyle\phi(s)\stackrel{\text{\rm def}}{=\!=} \lim\limits
_{\Delta \downarrow 0}\frac{\mathsf{P} \{\mathscr{T} \in
(s,s+\Delta]| \mathscr{T}>s\}}{\Delta} =
\frac{\Phi'(s+0)}{1-\Phi(s)}
$$
is called {\it the intensity} of the end of the random period
$\mathscr{T}$ at the time $s$ given $\mathscr{T}\geqslant s$.
\hfill \ensuremath{\triangleright}
\end{definition}

\begin{remark}
Definition \ref{zv1} holds:

$
\mathsf{P} \{\mathscr{T}\in(s,s+\Delta]| \mathscr{T}>s\}=\phi(s)\Delta+ o(\Delta).
$\hfill \ensuremath{\triangleright}
\end{remark}
\begin{remark}\label{zv9} 
C.d.f. and its density can be reconstructed by the intensity:
\begin{equation}\label{zv2} 
\Phi(s)=1-\exp\left(-\int\limits_0^s \phi(u)\,\mathrm{d}\, u\right);
\end{equation}
\begin{equation}\label{zv3} 
\Phi'(s)= \phi(s)\exp\left(-\int\limits_0^s \phi(u)\,\mathrm{d}\, u\right).
\end{equation}

It is easy to see: if c.d.f. is exponential then its intensity is a constant.

If c.d.f. $\Phi(s)$ is not continuous, then it can also use the intensity of c.d.f. $\Phi(s)$, but the formulae (\ref{zv2}) and
(\ref{zv3}) are not true. \hfill \ensuremath{\triangleright}
\end{remark}

So, the behaviour of considered queueing system with infinitely many servers, will be described on the intensities' language:
\\
{\bf 1.} The intensity of the input of next customer at the time moment $t$ is $\lambda(t)$, i.e..
$$
\mathsf{P} \left\{
\begin{array}{l}
\mbox{Next customer comes}
\\
\mbox{in the time interval $(t, t+\Delta]$}
\end{array}
\right \}=\lambda(t)\Delta+ o(\Delta).
$$
\\
{\bf 2.} Denote $h_i(t)$ -- the service intensity (at the time $t$) of $i$-th (by the order of arrivals) customers given at the time $t$ in the queueing system there are $n$ customers, i.e.
$$
\mathsf{P} \left\{
\begin{array}{l}
\mbox{$i$-th customer being on the service}
\\
\mbox{at the time $t$ ends own service}
\\
\mbox{in the time interval $(t,t+\Delta]$ }
\end{array}
\right\}
=h_i(t)\Delta+ o(\Delta).
$$

{\it Here and everywhere the $i$-th customer being in the queueing system, is the $i$-th by order of arrival from all customers being in the queueing system.}

In the case when $\lambda(t)\equiv \lambda$ and $h_i(t)\equiv \mu$, the considered queueing system is a queueing system $M|M|\infty$.
Firstly, this queueing system was studied in \cite{Zv8}; therein the following statement has been proved:
if $\mu<\lambda$, then the distribution of the quantity $n_t$ of the customers being in the queueing system converges to the stationary distribution:

\hspace{2cm}$
\displaystyle\lim\limits _{t\to\infty}\mathsf{P} \{n_t=k\}=\mathscr{P} _k=\mathscr{P} _0\frac{\lambda^k}{\mu^k}, \mbox{ where } \mathscr{P} _0=1-\frac{\lambda}{\mu}.
$

Later, this result was generalized for the case when the {\it intensity} of the input flow depends on the quantity $n_t$ of the customers being in the queueing system, i.e.
$\lambda(t)\equiv \lambda_k$ for $n_t=k$ (see, e.g., \cite[Ch.4,
\S 4-5]{zv4}):
\begin{equation}\label{zv4} 
\lim\limits _{t\to\infty}\mathsf{P} \{n_t=k\}=\mathscr{P} _k=\displaystyle \mathscr{P} _0\frac{\displaystyle \prod\limits_ {i=0} ^{k-1}\lambda_i}{\mu^k}, \mbox{ where } \mathscr{P} _0=\left(\sum\limits _{k=0}^\infty \frac{\displaystyle \prod\limits_ {i=0} ^{k-1}\lambda_i} {\mu^k}\right)^{-1}.
\end{equation}
Moreover, if the intensity of the service of all customers being in the queueing system depends on the quantity $n_t$ of the customers being in the queueing system (i.e. on $n_t$ ), i.e. $h_i(t)
\equiv \mu_k>0$ for $n_t=k$, then the formula (\ref{zv4}) is true, but with replacing of $\mu^k$ by $\prod\limits_ {i=1} ^{k}\mu_i$.
Cleary, if $\lambda_k=0$ for some $k$, then (\ref{zv4}) is useful for the queueing system with bounded quantity of servers, without places for the queue (without the buffer to wait for service).

In the case described above, when
\begin{equation}\label{zv5} 
\begin{array}{l}
\lambda(t)\not\equiv \lambda,\qquad \mu(t)\not\equiv \mu,
\\
\lambda(t)=\lambda_{n_t}, \qquad h_i(t)=\mu_{n_t},
\end{array}
\end{equation}
where $n_t$ is a quantity of the customers being in the queueing system at the time moment $t$, the (random) service times, and the input flow {\it are not independent}.

In \cite{Zv11} the queueing system with Poisson input and arbitrary absolutely continuous distribution of time of service (queueing system $M|G|\infty$) was studied. Here all customers have the same distribution of the service. The length of service and input flow are independent. In the language of the intensities, it means that if at the time moment $t$ in the system there are $n_t$ customers, and elapsed time of the service of $i$-th is $x^{(i)}_t$ ($i=1,2,\ldots,n$), then the intensity of service of $i$-th customer is equal $h_i(t)=h\left(x^{(i)}_t\right)$ at the time moment $t$.
.

Denote $\mu\stackrel{\text{\rm def}}{=\!=} \left(\displaystyle \int\limits_0^\infty s\,\mathrm{d}\, G(s)\right)^{-1}$, where $G(s)=1-\exp\left(-\displaystyle \int\limits_0^s h(u)\,\mathrm{d}\, u\right)$ -- c.d.f. of the service length. 
 The statement was made (see \cite{Zv11}): for $k>0$,
\begin{multline*}
\lim\limits _{t\to\infty} \mathsf{P} \{n_t=k; x_t^{(1)}<\alpha_1,
x_t^{(2)}<\alpha_2,\ldots x_t^{(k)}<\alpha_k \}=
\\
=\mathscr{P} (k;\alpha_1, \alpha_2, \ldots, \alpha_k)=\mathscr{P} _0\prod\limits_{i=1}^k \lambda_i G(1-\alpha_i),
\end{multline*}
where the normalizing multiplier $\mathscr{P} _0$ is calculated in the same way as in (\ref{zv4}).

This fact for the queueing system $M|G|m|0$, but in a more generalized form, without assumption about the absolutely continuity of c.d.f. of the service time $G(s)$, was proved in \cite{seva1} (see also \cite{seva2}).

So,  the studied queueing system had (see \cite{Zv11})
\begin{equation}\label{zv6} 
\lambda(t)\equiv \lambda, \qquad h_i(t)=h_i\left(x_t^{(i)}\right),
\end{equation}
where $x_t^{(i)}$ is the elapsed time of service of the $i$-th by the order of the arrival customer ($i=1,2,\ldots,n_t$, where at the time $t$, the queueing system contains $n_t$ customers).

Later, the results based on \cite{Zv8}, was obtained in many papers (see, e.g., \cite{zv2,Zv9,Zv12,bor}).

The natural generalization of the cases of (\ref{zv5}) and (\ref{zv6})
is the queueing system in which the intensity of input flow and the intensities of service depend not only of $n_t$, but also they depend of the elapsed service times and of the time from the last customer arrival:

\hspace{0.3cm}$
\lambda(t)=\lambda\left(x_t^{(0)}; x_t^{(1)}, x_t^{(2)},\ldots,x_t^{(n_t)}\right), \;\; h_i(t)=h_i\left(x_t^{(0)}; x_t^{(1)}, x_t^{(2)},\ldots,x_t^{(n_t)}\right),
$
\\
where $x_t^{(0)}$ is the time from the last arrival of the customer (at the time $t$), and $x_t^{(i)}$ is the elapsed service time of $i$-th by the order of arrival customer staying in the queueing system at the time $t$.
The knowing of the stationary distribution of the queueing system is very important for applications of the queueing theory; this distribution is the base for calculation of many important parameters of the complex technical systems.

Also for applications, it is very important to know the rate of convergence of the parameters of some technical systems to the stationary value.
For this aim, it is enough to know the convergence rate of the distribution of the queueing system to the stationary one.

Such problems was considered, e.g., in \cite{ver1,ver2,zv7}.

For example, in \cite{zv7} the following fact was proved: for queueing system
with one server and infinitely many places in queue, in the case, when the intensity of the service $h$ and intensity of the input flow $\lambda$ depend on the quantity of customers in the queueing system $n_t$, and on the elapsed service time $x_t$, the next Theorem \ref{zv7} is proved.
\begin{theorem}\label{zv7} 
Denote $X_t\stackrel{\text{\rm def}}{=\!=} (n_t,x_t)$, where $n_t$ is
the quantity of customers in the queueing system at the time moment
$t$, $x_t$ is the elapsed service time (of the customer being on the service) at the time moment $t$.
In the case $n_t=0$, the intensity of the input flow is constant, and it is equal to $\lambda_0>0$, and $X_t\stackrel{\text{\rm def}}{=\!=}(0,0)$.

If $\lambda(n_t,x_t)$ and $h(n_t,x_t)$ are Borel measurable and bounded, and
$$
\begin{array}{l}
\sup\limits _{n,x:\;n>0}\lambda(n,x)= \Lambda<\infty;\qquad
\inf\limits _{n>0, x\geqslant 0} h(n,x)\geqslant\displaystyle
\frac{C_0}{1+x};
\\ \\
C_0>4(1+2\Lambda), \mbox{ and for some } k>1 \quad C_0>2^{k+1}(1+\Lambda 2^k),
\end{array}
$$
then the distribution of the process $X_t=(n_t,x_t)$ converges to the stationary distribution $\mathscr{P} $ on the state space $\mathscr{X}
= \{(0,0)\} \bigcup \{(n,x): \; n \in \mathbb N, \; x\ge 0\}$, and there exists $m>k$ and $C>0$, such that fo all $t\geqslant 0$
$$
\sup\limits _{S\in \mathscr{B} (\mathscr{X} )}|\mathsf{P} \{X_t\in S|X_0=(n,x)\}-\mathscr{P} (S)|\leqslant C\,\frac{(1+n+x)^m}{(1+x)^k},
$$
where $\mathscr{P} $ is a stationary
distribution of the process $X_t$.
\hfill \ensuremath{\triangleright}
\end{theorem}
In the poof of this Theorem \ref{zv7} there is (not full) algorithm of the calculation of the constant $C$.

But there exists many results about the {\it type} of the convergence rate of the distribution of the state of the queueing system to the stationary distribution (see, e.g., \cite{zv1, Zv14,bor, GK} et al.).
The main result of this works is following: some parameters of the queueing system $\chi_t$ converge to the stationary value $\tilde\chi$ with the rate $\phi(t)$, i.e. there exists the constant $K$ such, that for all $t$, $ |\chi_t-\tilde\chi|<K \phi(t) $ is true (here $\phi(t)\searrow 0$).

For example, for the queueing system studied in \cite{zv7}
(Theorem \ref{zv7}), for the convergence of the distribution of the process $X_t$ to the  stationary distribution with the rate less then $\displaystyle \frac{K}{(1+x)^k}$, $k\in[0,C_0-1)$, it is sufficiently $\Lambda>C_0$, but without additional conditions of the Theorem \ref{zv7} there was no bounds for the constant $K$.

Emphasize, the queueing system $MI|GI|\infty$ is well studied, and many its parameters are well known (see, e.g., \cite{Zv9}). For example, this is $\mathscr{P}_k(t)$ (the probability that the system contains $k$ customers in the service at the time moment $t$), the expectation and the variance of the busy period, and the Laplace transform of c.d.f. of the busy period (it at the time $t=0$ the queueing system $MI|GI|\infty$ is idle).

Correspondingly, these results make it easy to estimate the convergence rate of $\mathscr{P} _k(t)$ to their stationary values $\mathscr{P} _k$ -- see, for example, \cite{Zv3}.

\subsection{Well-known facts and the aim of this paper}
Recall that the queueing system $MI|GI|\infty$ consists of the infinite number of servers operating independently, and the random time of service of all customers has the same distribution with c.d.f. $G(s)$.
Input flow is Poisson flow with constant parameter $\Lambda$.
At the initial time moment $t=0$ the queueing system is idle (no customer on the service).

Our goal is to study the behaviour of the queueing system with infinitely many servers, in which the input flow and (random) service times are described by the intensities depending on the {\it  queueing system full state}.
Therefore, the service times and input flow are dependent.
Naturally, despite this dependence, it is assumed that the characteristics of r.v. ``service times'' and ``the time interval between the arrivals of the customers'' are in some boundaries.
For example, the requirement that the incoming flow intensity does not be zero and bounded is naturally.

In this paper, the {\it  queueing system full state} at the time $t$ consists of time from the last arrival of the customer before time $t$, and the elapsed times of the service of the customer being on the service at the time $t$.

Suggested above dependence of the service times and input flow allows us to claim that the process of service  such queueing systems is Markov and {\it regenerative}.

The regeneration points of this process are the time moments when all components of the {\it  queueing system full state} are equal to zero.
It is moments when in the idle queueing system the new customer comes.
I.e. the regeneration point is the time moment when the queueing system serves only one customer, and the elapsed time of its service is equal zero; also at this time moment, the time from the last arrival is equal zero.

\section{Description of the studied queuing system }\label{zv16} 
So, consider the queueing system, in which the intensity of input flow and intensities of service depend on the {\it  queueing system full state}, which will be described below.

This {\it  queueing system full state} at the time moment $t$ includes: the time $x^{(0)}_t$ from the last arrival of the customer.
And elapsed times of service $x^{(i)}_t$ of $i$-th customer from the customer being in the queueing system at the time $t$; for convenience, also the {\it  queueing system full state} includes the variable $n_t$ -- the quantity of the customers in the queueing system at the time moment $t$. ({\it The customers are numbered by the order of the busy periods arrivals.})
So, the {\it  queueing system full state} is a vector $X_t=\left(n_t, x^{(0)}_t;
x^{(1)}_t,x^{(2)}_t,\ldots,x^{(n_t)}_t\right)$.
As $i$-th customer arrived before $(i+1)$-th one, cleary, that $x^{(i)}_t\geqslant
x^{(i+1)}_t$ for $i=1,\ldots,n_t$.
Also $x^{(0)}_t\leqslant x^{(n_t)}_t$, i.e. $n_t$-th customer appeared in the system no later than the time $t-x^{(0)}_t$.

The intensity of the input flow is function from $X_t$:\; $\lambda=\lambda(X_t)$ -- as soon as the intensity of service of $i$-th customer is a function from $X_t$: $h_i=h_i(X_t)$.
Recall, the $i$-th customer has the elapsed service time $x^{(i)}_t$.
Hence, for little time interval $\Delta>0$ given $X_t=\left(n_t,
x^{(0)}_t; x^{(1)}_t,x^{(2)}_t,\ldots,x^{(n_t)}_t\right)$
(where $x^{(i)}_t\geqslant x^{(i+1)}_t$) the formulae (\ref{zv8}) are true:
\begin{equation}\label{zv8} 
\left.
\begin{array}{l}
\mathsf{P} \left\{
\begin{array}{l}
n_{t+\Delta}=n_t+1,\; x_{t+\Delta}^{(0)}=
x^{(n_{t+\Delta})}_{t+\Delta}\in(0;\Delta),
\\
\mbox{and $x_{t+\Delta}^{(i)}= x_{t}^{(i)}+\Delta$ for all $i=1,\ldots,n_t$ }
\\
x^{(n_t+1)}_{t+\Delta}\in (0;\Delta)
\end{array}
\right\} =\lambda(X_t)\Delta+o(\Delta);
\\
\mathsf{P} \{n_{t+\Delta}=n_t-1\}=\displaystyle \sum\limits
_{i=1}^{n_t}h_i (X_t)\Delta+o(\Delta);
\\
\mathsf{P} \left\{
\begin{array}{l}
n_{t+\Delta}=n_t-1, \mbox{ and}:
\\
\mbox{1. for all $j<i$, } x_{t+\Delta}^{(j)}=x_{t}^{(j)}+\Delta,
\\
\mbox{2. for all $j>i$ } x_{t+\Delta}^{(j)}=x_{t}^{(j+1)}+\Delta
\end{array}
\right\}=h_i(X_t)\Delta+o(\Delta)
\\
\hspace{5cm}\mbox{ -- for all $i=1,2,\ldots,n_t$};
\\
\mathsf{P} \left\{n_{t+\Delta}=n_t; x_{t+\Delta}^{(i)}=x_{t}^{(i)}+\Delta \mbox{ for all $i=0,\ldots,n_t$}\right\}=
\\
\hspace{4cm}=1-\left(\lambda(X_t)+\displaystyle \sum\limits
_{i=1}^{n_t}h_i (X_t)\right)\Delta+o(\Delta).
\end{array}
\right\}
\end{equation}

From (\ref{zv8}) it follows, that the distribution of $X_{t+\Delta}$ depends only on distribution of $X_t$ and $\Delta$. So, $X_t$ is Markov process on the state space $\displaystyle \mathscr{X} \stackrel{\text{\rm def}}{=\!=} \bigcup\limits_{i=0}^\infty\mathscr{S} _i, $ where $\mathscr{S} _i\stackrel{\text{\rm def}}{=\!=} \{i\}\times\prod\limits_{j=0}^i \mathbb{R}_+$, $i\in\mathbb{Z}_+$.
Here the set $\mathscr{S} _0$ is the set of the idle states of the queueing system (when there no customers in the queueing system).

Denote $\mathscr{P} _t^{X_0}$ -- the distribution the process $X_t$ with the initial state $X_0$: $\mathscr{P} _t^{x}(S)=\mathsf{P} \{X_t\in S|X_0=x\}$, $x\in \mathscr{X} $, $S\in\mathscr{B}( \mathscr{X})$.

Recall, that the time moments when the process $X_t$ hits to the state $(1,0;0)$, i.e. when the process $X_t$ exits the set $\mathscr{S} _0$ (in the idle queueing system new customer comes), are the regeneration times for the process $X_t$.
Thus, if the time between two regeneration times has the finite expectation, then the process $X_t$ is ergodic, and its distribution weakly converges to the some stationary (invariant) distribution.

\begin{remark}
The weak class of the queueing system may be described by this way.
For example, the queueing system $M|M|\infty$, $M|G|\infty$, $M|G|m|0$, and other.\hfill
\ensuremath{\triangleright}
\end{remark}

For the well-studied queueing system $M|G|\infty$, the finite expectation of the service time is enough for the finiteness of the regeneration period.
And in this case the convergence rate of the parameters $\mathscr{P} _k(t)$ to their stationary values $\mathscr{P} _k$ increases if the number of finite moments of the service time increases -- see, e.g., \cite{Zv3}.

Accordingly, it is naturally to suppose that the (dependent) service times have finite moments, and the parameters of the input flow can be bounded by the parameter of Poisson input flow.

Suppose the conditions:
\begin{equation}\label{zv10} 
\begin{array}{l}
h_i(X_t)\geqslant \displaystyle \frac{K}{1+x^{(i)}_t},\;\; K>2;\quad
0<\lambda_0\leqslant\lambda(X_t)\leqslant \Lambda<\infty.
\end{array}
\end{equation}

It holds from the first condition of (\ref{zv10}) that c.d.f. of the service time $G_i(s)$ satisfies the inequality $ \displaystyle G_i(s)\geqslant
1-\displaystyle \frac{1}{(1+s)^K}, $ i.e. the service time has
$k\in[0;K)$ finite moments (see Remark \ref{zv9}).

Also the density of c.d.f. of the intervals between arrivals of the customers
$F(s)$ satisfies the inequality $ \displaystyle f_i(s)\leqslant
\Lambda e^{-\lambda_0 s}, $ and distribution of the time intervals $\eta_i$ between arrivals of the customers has any moments: 
\begin{equation}\label{zzzv}
\displaystyle \mathbb{E} \,(\eta_i)^k\leqslant\displaystyle
\frac{k!\Lambda}{\lambda_0^{k+1}}, 
\end{equation}
in particular, this inequality is correct
for the {\it idle periods} of $X_t$, i.e. time intervals, where
$n_t\equiv 0$.\hfill \ensuremath{\triangleright}

Our goal is to find the computable bounds of the rate of convergence of the distribution of the process $X_t$ to the stationary distribution in conditions (\ref{zv10}).

The {\it coupling method}  will be used for this aim -- see Section
\ref{zv11}.

Previously, recall some information from the queueing theory, and give some necessary inequalities.

\section{Some information about the queueing sys\-tem\\$MI|GI|\infty$}\label{zv20} 
Here there are well-known results about the queueing system $MI|GI|\infty$ (see \cite{Zv9,Zv10} at al.).

The input flow of $MI|GI|\infty$ is Poisson, its intensity is constant $\Lambda$. The arrived $i$-th customer immediately became be serviced.
Its service time is r.v. $\xi_i$ with c.d.f. $G(s)$.
All r.v. (service times of and times between the customer arrivals) are independent mutually.

\textit{At the initial time $t=0$ the queueing system is idle.}

Suppose, that $ \displaystyle \mathbb{E}
\,(\xi_i)^k=k\int\limits_0^\infty s^{k-1}(1-G(s))\,\mathrm{d}\,
s=m_k<\infty $ for some $ K\geqslant 1$.
For example,
\begin{equation}\label{zv21} 
G(s)=1-\displaystyle \frac{1}{(1+s)^{K}},\qquad K>2.
\end{equation}

The condition (\ref{zv21}) guarantees that for all $k\in[0,K)$, $\mathbb{E}
\,(\xi_i)^k<\infty$.
It is obvious that
$
\displaystyle \mathbb{E} \,\xi_i=\frac{1}{K-1};\quad \mathbb{E} \,(\xi_i)^2=\frac{2}{K^2-3K+2};\mbox{ etc.}
$

The loading of the queueing system is $\rho\stackrel{\text{\rm
def}}{=\!=} \Lambda m_1$.

The behaviour of the process of the queueing system described by Markov process $\stackrel{\circ}{X}_t=\left(\stackrel{\circ}n_t, \stackrel{\circ}x^{(0)}_t; \stackrel{\circ}x^{(1)}_t,\stackrel{\circ}x^{(2)}_t,\ldots,\stackrel{\circ}x^{(\stackrel{\circ}n_t)}_t\right)$ anew , here $\stackrel{\circ}n_t$ is a quantity of the customers in the queueing system at the time moment $t$, and $\stackrel{\circ}x^{(0)}_t$ is the time from a last arrival of the customers to the queueing system; $\stackrel{\circ}x^{(i)}_t$ is the elapsed time of the service of $i$-th customer being in the queueing system ($i=1,2,\ldots, \stackrel{\circ}n_t$).
Here and everywhere else the sign $\circ$ over a random variable, process or numerical characteristics indicates the ``classic'' queueing system $MI|GI|\infty$.)

Denote $\stackrel{\circ}{\mathscr{P}}_t$ -- the distribution the process $\stackrel{\circ}{X}_t$: $ \stackrel{\circ}{\mathscr{P}}
_t(S)=\mathsf{P} \{\stackrel{\circ}{X}_t\in S\} $ for all
$S\in\mathscr{B} (\mathscr{X} )$, and again $\mathscr{X}
\stackrel{\text{\rm def}}{=\!=}
\bigcup\limits_{i=0}^\infty\mathscr{S} _i, $ where $\mathscr{S}
_i\stackrel{\text{\rm def}}{=\!=} \{i\}\times \prod\limits_{j=0}^i
\mathbb{R}_+$, $i\in\mathbb{Z}_+$.

In condition (\ref{zv21}), it is known (see, e.g., \cite{Zv9}), that 
$$
\displaystyle\mathcal P_k(t)\stackrel{\text{\rm def}}{=\!=}
\mathsf{P} \{\stackrel{\circ}n_t=k\}=\stackrel{\circ}{\mathscr{P}}
_t(\mathscr{S} _0)=\displaystyle \frac{e^{-\mathfrak{G}
(t)}\big(\mathfrak{G} (t)\big)^k}{k!}, $$ 
where $\displaystyle\mathfrak{G} (t)\stackrel{\text{\rm def}}{=\!=}
\displaystyle \Lambda\int\limits_0^t(1-G(s))\,\mathrm{d}\, s. $
\begin{remark}\label{zv26} 
$\mathfrak{G} (+\infty)=\Lambda m_1$, and $\mathcal P_0(t)\geqslant e^{-\Lambda m_1}=e^{-\rho}$ for all $t\geqslant 0$.\hfill \ensuremath{\triangleright}
\end{remark}

Moreover, the distribution of the busy period of the queueing system $M|G|\infty $ is known (see, e.g., \cite{Zv10}).
Denote $\stackrel{\circ}\zeta_i$ -- $i$-th busy period of the queueing system
$M|G|\infty $.
This r.v. has c.d.f. $ \displaystyle
B(x)\stackrel{\text{\rm def}}{=\!=} \mathsf{P}
\{\stackrel{\circ}\zeta_i\leqslant x\}=1-
\frac{1}{\Lambda}\displaystyle \sum\limits _{k=1}^\infty
c^{n\ast}(x), $ where $c(x)\stackrel{\text{\rm def}}{=\!=}
\displaystyle \Lambda(1-G(x))e^{-\mathfrak{G} (x)}$, and $c^{n\ast}$
is $n$-th convolution of the function $c(x)$.

Also the Laplace transform of the distribution $B(x)$ is known (see, e.g., \cite{Zv10}):
\begin{equation}\label{zv17} 
\mathcal{L}[B](s)=1+\frac{s}{\Lambda}-\frac{1}{\Lambda \displaystyle \int\limits_0^\infty \exp\left(-st-\Lambda\int\limits_0^t [1-G(v)]\,\mathrm{d}\, v\right)\,\mathrm{d}\, t}.
\end{equation}
Moreover, there exists the formulae for the moments of the busy period:
\begin{equation}\label{zv18} 
\mathbb{E}
\,(\stackrel{\circ}\zeta_i)^n=(-1)^{n+1}\left\{\frac{e^\rho}{\Lambda}nC^{(n-1)}
- e^\rho \displaystyle \sum\limits _{k=1}^{n-1}C_n^k \mathbb{E}
\,(\stackrel{\circ}\zeta_i)^{n-k} C^{(k)} \right\}, \qquad
n\in\mathbb{N},
\end{equation}
where
$\displaystyle
C^{(k)}=\int\limits_0^\infty (-t)^k \left(\exp\left(-\Lambda
\int\limits_0^\infty [1-G_0(v)]\,\mathrm{d}\, v\right)\right)\Lambda
[1-G_0(t)]\,\mathrm{d}\, t.
$

From (\ref{zv17}) and (\ref{zv18}), it can obtain
\begin{equation}\label{zv36} 
\begin{array}{l}
\hspace{2cm}\displaystyle \mathbb{E} \, \stackrel{\circ}\zeta_i=
\frac{e^\rho-1}{\Lambda} \mbox{ and }
\mathbb{E} \,(\stackrel{\circ}\zeta_i)^2=
\\ =\displaystyle \frac{2e^ {2\rho}}{\Lambda}\displaystyle \int\limits_0^\infty \left(\exp\left(-\Lambda \int\limits_0^\infty [1-G(v)]\,\mathrm{d}\, v\right)\right)\,\mathrm{d}\, t+ \displaystyle \frac{2e^\rho} {\Lambda^2}-\frac{2}{\Lambda^2}-2\displaystyle \frac{e^\rho-1}{\Lambda^2}.
\end{array}
\end{equation}

But for calculation of $\mathbb{E} \,(\stackrel{\circ}\zeta_i)^k $, $k>2$, the use of the formulae (\ref{zv17}) and (\ref{zv18}) is very complicated and hardly to apply in practice.

But for our aim, we need to use the bounds for $\mathbb{E} \,(\stackrel{\circ}\zeta_i)^k $.
For this goal, we use
\begin{proposition}\label{zv19} 
If $\mathbb{E} \,\xi_i^k<\infty$ for some $k\in\mathbb{N} $, then the follows ibequality is true:
\begin{equation}\label{zv25} 
\mathbb{E} \,(\stackrel{\circ}\zeta_i)^k\leqslant \frac{\mathbb{E} \,\xi_i ^k}{1-e^{-\rho}} \times \varphi((1-e^{-\rho}),k-1),
\end{equation}
where $\displaystyle \varphi(x,k)\stackrel{\text{\rm def}}{=\!=}
\displaystyle \sum\limits _{n=1}^\infty
n^{k-1}x^n=\left(x\,\frac{\,\mathrm{d}\,}{\,\mathrm{d}\,
x}\right)^k\frac{1}{1-x}$. \hfill \ensuremath{\triangleright}
\end{proposition}

\noindent\textbf{Proof.} The function $c(x)$ from (\ref{zv17}) is nonnegative, and
\begin{multline*}
\displaystyle \int\limits_0 ^\infty c(x)\,\mathrm{d}\,
x=\int\limits_0 ^\infty \Lambda(1-G(x))\exp\left(\displaystyle
-\Lambda\int\limits_0^x(1-G(s))\,\mathrm{d}\, s\right)
\,\mathrm{d}\, x =
\\
=\left(1-\exp\left(\displaystyle
-\Lambda\int\limits_0^\infty(1-G(s))\,\mathrm{d}\, s\right)\right)=
\left(1-e^{\displaystyle -\Lambda
m_1}\right)=(1-e^{-\rho})\stackrel{\text{\rm def}}{=\!=} \varrho,
\end{multline*}
where again $\rho=\displaystyle \Lambda m_1$.

I.e. the function $\varsigma(x)\stackrel{\text{\rm def}}{=\!=}
\varrho^{-1}\, c(x)$ ia a density of c.d.f. of some r.v. $\vartheta$, and $\varsigma^{n\ast}(x)$ is a density of c.d.f. of
$\displaystyle \displaystyle \sum\limits _{j=1}^n\vartheta_j$ --
the summ of $n$ i.i.d. r.v. with the density of c.d.f. $c(x)$.

Thus, $\displaystyle \int\limits_0^ \infty
\varsigma^{n\ast}(x)\,\mathrm{d}\, x=1$.
Therefore,
\begin{multline*}
\mathbb{E} \,\stackrel{\circ}\zeta_i= \int\limits_0^\infty (1-B(x))\,\mathrm{d}\, x=\displaystyle \sum\limits _{n=1}^\infty \frac{1}{\Lambda}\int\limits_0^\infty c^{n\ast}(x)\,\mathrm{d}\, x =
\\
= \frac{1}{\Lambda}\displaystyle \sum\limits _{n=1}^\infty \varrho^n\int\limits_0^\infty \varsigma^{n\ast}(x)\,\mathrm{d}\, x=\frac{1}{\Lambda}\displaystyle \sum\limits _{n=1}^\infty \varrho^n =\frac{1}{\Lambda}\times \frac{\varrho}{1-\varrho}=\frac{e^{\rho}-1}{\Lambda}.
\end{multline*}

Now, if $m_{k+1}=\mathbb{E} \,\xi_i^{k+1}<\infty$ for some
$k\in\mathbb{N} $, then
\begin{multline*}
\mathbb{E} \,\vartheta^k=\int\limits_0 ^\infty x^k\varsigma(x)\,\mathrm{d}\, x =\int\limits_0 ^\infty x^k \varrho^{-1} \Lambda(1-G(x))e^{-\mathfrak{G} (x)}\,\mathrm{d}\, x \leqslant
\\
\leqslant \frac{\Lambda}{\varrho }\int\limits_0^ \infty x^k (1-G(x))\,\mathrm{d}\, x =\frac{\Lambda\mathbb{E} \,\xi_i^{k+1}}{\varrho(k+1)},
\end{multline*}
and $ \displaystyle \int\limits_0 ^\infty
x^k\varsigma^{n\ast}(x)\,\mathrm{d}\, x=\mathbb{E}
\left(\displaystyle \sum\limits
_{j=1}^n\vartheta_j\right)^k\leqslant n^{k}\mathbb{E}
\,\vartheta^k\leqslant \frac{ n^{k}\Lambda m_{k+1}}{\varrho (k+1)}.
$

In the last inequality the Jensen's inequality was used in the form: for
$k\geqslant 1$,
\begin{equation}\label{zv33} 
\displaystyle (a_1+\ldots+a_n)^k\leqslant n^{k-1}(a_1^k+\ldots +a_n^k).
\end{equation}

Now,
\begin{multline*}
\mathbb{E} \,(\stackrel{\circ}\zeta_i)^k=\int\limits_0 ^ \infty k
x^{k-1}(1-B(x))\,\mathrm{d}\, x =\frac{1}{\Lambda}\displaystyle
\sum\limits _{n=1}^\infty\int\limits_0^ \infty k x^{k-1}
c^{n\ast}(x)\,\mathrm{d}\, x =
\\
= \frac{1}{\Lambda}\displaystyle \sum\limits
_{n=1}^\infty\int\limits_0^ \infty k x^{k-1}
\varrho^n\varsigma^{n\ast}(x)\,\mathrm{d}\, x=
\frac{1}{\Lambda}\displaystyle \sum\limits _{n=1}^\infty \varrho^n
k\, \mathbb{E} \left(\displaystyle \sum\limits _{j=1}^n
\vartheta_j\right)^ {k-1}\leqslant
\\
\leqslant \frac{1}{ \Lambda} \displaystyle \sum\limits _{n=1}^\infty
\varrho^nk n^{k-1}\frac{\Lambda\mathbb{E} \, \xi_i^k}{\varrho k}=
\frac{\mathbb{E} \,\xi_i ^k}{\varrho } \displaystyle \sum\limits
_{n=1}^\infty n^{k-1} \varrho^n= \frac{\mathbb{E} \,\xi_i
^k}{\varrho } \times \varphi(\varrho,k-1).
\end{multline*}
Proposition \ref{zv19} is proved.
\hfill \ensuremath{\blacksquare}

\section{Bounds for the process $X_t$}\label{zv31} 
Return to the studied process $X_t$, described in the Section \ref{zv16}.
Suppose, that at the initial time $t=0$, $X_0=(0,0)$.
I.e. the queueing system is idle, and it waits the first customer.
At the time of arrival of the first customer, the process $X_t$ hits to the state $(1,0;0)$; and this time moment is the regeneration time.

Recall, that the process $\stackrel{\circ}{X}_t$ described in the Section \ref{zv20} also starts from the state $(0,0)$: $\stackrel{\circ}X_0=(0,0)$.
At the time of arrival of the first customer, the first regeneration period begins.
\begin{definition}\label{zzv}
The random variable $\eta$ does not exceed the random variable $\theta$
\textit{by distribution}, if for all $s\in\mathbb{R}$ the
inequality $ \displaystyle F_\eta(s)=\mathsf{P} \{\eta\leqslant s\}
\geqslant \mathsf{P} \{\theta\leqslant s\}=F_\theta(s) $ is true.
Hence, c.d.f. of r.v. $\theta$ does not exceed c.d.f. of r.v. $\eta$.
This is order relation. Denote it by $\eta\prec
\theta$.\hfill \ensuremath{\triangleright}
\end{definition}
\begin{proposition}\label{zv22} 
$\eta\prec \theta$ iff on the some probability space $(\Omega,\mathscr{F} ,\mathsf{P} )$ there exists r.v. $\eta'$ and $\theta'$ such, that
$\eta'\leqslant\theta'$ for all $\omega\in\Omega$, and $\mathsf{P}
\{\eta\leqslant s\}=\mathsf{P} \{\eta'\leqslant s\}$ and $\mathsf{P}
\{\theta\leqslant s\}=\mathsf{P} \{\theta'\leqslant s\}$ for all
$s\in\mathbb{R}$.\hfill \ensuremath{\triangleright}
\end{proposition}

\noindent\textbf{Proof.}
see \cite{sht}.
\hfill \ensuremath{\blacksquare}
\begin{proposition}\label{zv23} 
If $X_0=\stackrel{\circ}X_0=(0,0)$, and the condition
(\ref{zv10}) holds true (recall, that input flow for $\stackrel{\circ}{X}_t$
is Poisson with the intensity $\Lambda$, and the time of service
has c.d.f. (\ref{zv21})), then for all $t\geqslant 0$ condition (\ref{zv24}) is true:
\begin{equation}\label{zv24} 
n_t\prec \stackrel{\circ}n_t. 
\end{equation}
\hfill \ensuremath{\triangleright}
\end{proposition}

\noindent\textbf{Proof.} The proof is based on the
\textit{method of one probability space}.
Namely, on the one probability space, the processes $X_t$ and $\stackrel{\circ}X_t$ can be defined (more specifically, their version) such, that for all trajectories of both processes, the inequality $n_t\leqslant \stackrel{\circ}n_t$ is true -- for all
$t\geqslant 0$.

This construction is possible by two reasons.

{\bf I.} The input flow of the process $X_t$ has the intensity, which does not exceed the intensity of the input flow of the process $\stackrel{\circ}X_t$
($\lambda(X_t)\leqslant \Lambda$ from (\ref{zv10})).

Therefore, on the one probability space,  two input flows for both process can be defined: the input flow $\mathfrak{F}_t$ with intensity $\lambda(X_t)$ for the process $X_t$, and the input flow $\mathfrak{F}^+_t$ with intensity $\Lambda-\lambda(X_t)$; the union of the flows $\mathfrak{F}_t$ and $\mathfrak{F}^+_t$ is the flow with intensity $\Lambda$ for the process $\stackrel{\circ}X_t$.

So, at any time $t\geqslant 0$, the quantity of the customers arriving to the queueing system with the process $X_t$, does not exceed the quantity of the customers arriving to the queueing system with the process $\stackrel{\circ}X_t$. Moreover, at the times of arrivals of the customers to the queueing system with the process $X_t$, contemporaneously the customers arrive into the queueing system with the process $\stackrel{\circ}X_t$; this customers are called ``common''.

In other words, at the any time moment, the vector $\left( x^{(1)}_t,x^{(2)}_t,\ldots,x^{(n_t)}_t\right)$ consists of the components of the vector $\left( \stackrel{\circ}x^{(1)}_t,\stackrel{\circ}x^{(2)}_t,\ldots,\stackrel{\circ}x^{(\stackrel{\circ}n_t)}_t\right)$.
I.e. the system described by the process $X_t$, the vector $\left(x^{(1)}_t,x^{(2)}_t,\ldots,x^{(n_t)}_t\right)$ consists of the elapsed service times of the customers. And these elapsed service times coincide with some components of the vector $\left( \stackrel{\circ}x^{(1)}_t,\stackrel {\circ}x^{(2)}_t,\ldots,\stackrel {\circ}x^{(\stackrel{\circ}n_t)}_t \right)$.

{\bf II.} The intensity of service of $i$-th customer being in the queueing system described by the process
$\stackrel{\circ}X_t=\left(\stackrel{\circ}n_t,
\stackrel{\circ}x^{(0)}_t; \stackrel{ \circ}x^{(1)}_t,
\stackrel{\circ}x^{(2)}_t,\ldots, \stackrel{\circ}x^{(\stackrel{
\circ}n_t)}_t\right)$, is equal $ \displaystyle
\stackrel{\circ}h_i(\stackrel {\circ}X_t)= \frac{G'\left(\stackrel{
\circ}x^{(i)}_t\right)} {1-G\left(\stackrel{ \circ}x^{(i)}_t\right)}
=\frac{K}{1+\stackrel{ \circ}x^{(i)}_t}. $
Thereby, the intensity of the service of the ``common'' customers in the process $\stackrel{\circ}X_t$ does not exceed the intensity of the service of the same customers in the process $X_t$.
This is the second reason to apply the method of one probability space.

Considering that the time of service has c.d.f. which can be calculated by the intensity of the service (see (\ref{zv2})--(\ref{zv3})), it can obtain the next fact.
The residual time of service of the ``common'' customers in the process $X_t$ does not exceed \textit{ by the distribution } the residual time of service of these customers in the process $\stackrel{\circ}X_t$.
So, it can construct on the one probability space this r.v. by such a way that the ``common'' customers will complete the service in the process $X_t$ no later than in the process $\stackrel{\circ}X_t$ -- see Proposition \ref{zv22}.

Therefore, in the queueing system described by the process $\stackrel{\circ}X_t$ the flow of the customers is larger than in the queueing system described by the process $X_t$, and simultaneously arrived (common) customers earlier leave the queueing system described by the process $X_t$.
This is the reason for the truth of the Proposition \ref{zv23}. \hfill \ensuremath{\blacksquare}
\begin{corollary}
From (\ref{zv24}), it follows that the property of orderliness in the distribution (see Definition \ref{zzv}) not true for only $n_t$ and $\stackrel{\circ}n_t$, but it is true also and for the busy periods of the processes $X_t$ and $\stackrel{\circ}X_t$.

The processes $X_t$ and $\stackrel{\circ}X_t$ are regenerative (recall that the regeneration points are the time points when the process is equal to $(1,0;0)$, i.e. the beginning of the busy period).
Hence, the distribution of the busy periods $\zeta_i$ of the process $X_t$ and the busy periods $\stackrel{\circ}\zeta_i$ of the process $\stackrel{\circ}X_t$ not depend on their number $i$.
So, $ \zeta_i\prec\stackrel{\circ}\zeta_t, $ and for all $k>0$ the inequality $ \displaystyle\boxed{\mathbb{E}
(\zeta_i)^k\leqslant\mathbb{E} (\stackrel{\circ}\zeta_i)^k} $ is true.
It means that the bounds (\ref{zv25}) are true for $\zeta_i$.\hfill
\ensuremath{\triangleright}
\end{corollary}
\begin{corollary}\label{zero} 
From Remark \ref{zv26}, and (\ref{zv24}) it follows that for all
$t\geqslant 0$ (and $X_0=(0,0)$) the inequality
\begin{equation}\label{zv32}
\mathsf{P} \{X_t\in\mathscr{S} _0\}=\mathsf{P} \{n_t=0\}\geqslant \mathcal{P}_0(t)\geqslant e^{-\rho}.
\end{equation}
is true.
\hfill \ensuremath{\triangleright}
\end{corollary}

\section{Coupling method}\label{zv11} 

The original coupling method for the Markov processes have been applied to the Markov chains with a finite number of states (see \cite{Zv6,Zv14}).
The main idea of this method is as follows.

Let two homogeneous Markov processes $Y_t$ and $Z_t$ (on the probability space $(\Omega,\mathscr{F} ,\mathsf{P} )$ with the filtration $\mathscr{F} _t$) with the same transition function, have the different states at the initial time moment.

Denote $\tau=\tau(Y_0,Z_0)\stackrel{\text{\rm def}}{=\!=}
\inf\{t>0:\, Y_t=Z_t\}$.
So, $Y_\tau=Z_\tau$, and for all $t>\tau$, the distributions of the processes $Y_t$ and $Z_t$ coincide in accordance with the Markov property, i.e. $ \mathsf{P} \{Y_t\in S|t>\tau\}=\mathsf{P} \{Z_t\in S|t>\tau\} $ for all
$S\in\mathscr{B} (\mathscr{X} )$, where $\mathscr{B} (\mathscr{X} )$ is $\sigma$-algebra on the state space $\mathscr{X}$ of $Y_t$ and $Z_t$.

{\it \bf The basic coupling inequality}:
\begin{multline}\label{zv12} 
|\mathsf{P} \{Y_t\in S\} - \mathsf{P} \{Z_t\in S\}|=
\\
=|\mathsf{P} \{Y_t\in S\;\&\; \tau>t\} - \mathsf{P} \{Z_t\in S\;\&\; \tau>t\} +
\\
+ |\mathsf{P} \{Y_t\in S\;\&\; \tau<t\} - \mathsf{P} \{Z_t\in S\;\&\; \tau<t\}|=
\\
=|\mathsf{P} \{Y_t\in S\;\&\; \tau>t\} - \mathsf{P} \{Z_t\in S\;\&\; \tau>t\} |\leqslant \mathsf{P} \{\tau>t\},
\end{multline}
as if $\tau<t $, then ``coupling'' of the processes $Y_t$ and $Z_t$ already occurred, and their distribution at the time $t$ are the same, and $\mathsf{P} \{Y_t\in S| \tau<t\}=\mathsf{P} \{Z_t\in S| \tau<t\}$; here $S\in\mathscr{B} (\mathscr{X} )$.

Therefore if there exists some bounds $\mathbb{E} \,
\phi(\tau(Y_0,Z_0))\leqslant C(X_0,Y_0)$ for some increasing positive function $\phi(t)$, then the Markov inequality is applicable:
\begin{equation}\label{zv13} 
\mathsf{P} \{\tau(Y_0,Z_0)>t\}=\mathsf{P} \{\phi(\tau(Y_0,Z_0))>\phi(t)\}\leqslant \frac{\mathbb{E} \, \phi(\tau(Y_0,Z_0))}{\phi(t)}\leqslant \frac{C(X_0,Y_0)}{\phi(t)}.
\end{equation}

Moreover, if there exists the stationary invariant distribution $\mathscr{P} $, and the distributions of $Y_t$ and $Z_t$ weakly converges to $\mathscr{P} $ as $t\to\infty$, then from (\ref{zv12}) and (\ref{zv13}) holds:
\begin{equation}\label{zv35} 
\left.\begin{array}{l} \hspace{1cm} |\mathsf{P} \{Y_t\in S\} -
\mathscr{P} (S)|\leqslant \displaystyle
\int\limits_\mathcal{\mathscr{X} }
\frac{C(X_0,y)}{\phi(t)}\mathscr{P} (\,\mathrm{d}\, y)=\displaystyle
\frac{\hat C(X_0)}{\phi(t)},
\\
\mbox{consequently,}
\\
\hspace{1cm} \|\mathscr{P} _t^{X_0}-\mathscr{P}
\|_{TV}\stackrel{\text{\rm def}}{=\!=} \sup\limits _{S\in
\mathscr{B} (\mathscr{X} )}|\mathscr{P} _t^{X_0}(S)-\mathscr{P}
(S)|\leqslant 2 \displaystyle \frac{\hat C(X_0)}{\phi(t)}.
\end{array}
\right\}
\end{equation}

Recall, that initially, the coupling method has been proposed for Markov chains with a discrete set of states (\cite{Zv6}), but studied in this paper process $X_t$ with transition low (\ref{zv8}) is the continuous process.
For two independent processes $X_t$ and $X_t'$ satisfying (\ref{zv8}), and with different initial states, $\mathsf{P} \{\tau(X_0,x_0')<\infty\}=0$.

For this reason, the notion of \emph{successful coupling} will be applied (see \cite{Zv5}).

\begin{definition} \label{zv30} 
Let $X_t$ and $X_t'$ are two independent versions of one homogeneous Markov process with the state space $\mathscr{X} $: $X_t$ starts from the state $X_0$, and $X_t'$ starts from the state $X_0'$.

\textsc {The successful coupling} (see \cite{Zv5}) of the processes $X_t$ and
$X_t'$ is the pair process $\mathcal{Z}_t=(Y_t,Y_t')$, defined on the some probability space, and this pair process satisfies the conditions:
\begin{itemize}
\item [(i)] for all $s\geqslant 0$ and $S\in \mathscr{B} (\mathscr{X} )$ the equalities

$
\mathbf{P} \{Y_t\in S\}= \mathbf{P} \{X_t\in S\}$, $ \mathbf{P} \{Y'_t\in S\}= \mathbf{P} \{X'_t\in S\} are true.
$
(Correspondingly, $Y_0=X_0$ and $Y'_0=X'_0$.)

{Emphasize, that for all fixed time moment $t\geqslant 0$ the random variables $Y_t$, $Y_t'$, $X_t$ and $X_t'$ are considered only as some random variables, not as the stochastic processes!
Here there is no coincidence of the finite-dimensional probability distributions of two processes, but only the coincidence of the marginal distributions.
}
\item[(ii)] for all $t>\tau=\tau(X_0,X'_0)=\tau(Y_0,Y'_0) \stackrel{\text{\rm def}}{=\!\!=\!\!=}\inf\{t\geqslant 0:\,Y_t=Y_t'\}$ the equality $Y_t=Y_t'$ is true.
\item[(iii)] For any initial states $X_0,X_0'\in \mathscr{X}$, $ \mathbf{P} \{\tau(X_0,X'_0)<\infty\}=1$.
\hfill \ensuremath{\triangleright}
\end{itemize}
\end{definition}
For \emph{successful coupling}, the \textit{ basic coupling inequality } (\ref{zv12}):
\begin{multline*}
|\mathsf{P} \{X_t\in S\}-\mathsf{P} \{X_t'\in S\}|=|\mathsf{P} \{Y_t\in S\}-\mathsf{P} \{Y_t'\in S\}|
\\
=|\mathsf{P} \{Y_t\in S\}-\mathsf{P} \{Y_t'\in S\}|\times (\mathbf{1} (\tau> t)+ \mathbf{1} (\tau\leqslant t))\leqslant \mathsf{P} \{\tau> t\}
\end{multline*}
is true.

So now, our goal is the construction of the successful coupling for two version of the process $X_t$ ($X_t$ with initial state $X_0=(0,0)$ and $X_t'$
with initial state $X_0'=\left(n_0', x_t'{}^{(0)};
x_t'{}^{(1)},x_t'{}^{(2)},\ldots,x_t'{}^{(n'_0)}\right)$), and both processes satisfy the conditions (\ref{zv8}).
From this construction of the successful coupling, the bounds for the convergence rate of the distribution of the process $X_t$ to the stationary distribution will be obtained.

For use the coupling method, the next Lemma is needful.
\subsection{Basic coupling Lemma}
This Lemma exists in many papers (see, e.g., \cite{Zv13,verbut}).
Here, this Lemma is in the simplest form.

\begin{definition} The common part of the distributions of two r.v. $\xi_1$ and $\xi_2$ having the c.d.f. $\Psi_j(s)=\mathbf{P}\{\xi_j\leqslant s\}$ $s\in \mathbb{R}$ ($j=1,2$) correspondingly, is
$\displaystyle \varkappa \stackrel{{\rm def}}{=\!\!\!=}
\varkappa(\Psi_1(s),\Psi_2(s)) \stackrel{{\rm def}}{=\!\!\!=}
\displaystyle \int\limits _{-\infty}^\infty
\min(\psi_1(u),\psi_2(u))\,\mathrm{d}\, u, $
where\\

$\psi_j(s) \stackrel{{\rm def}}{=\!\!\!=} \left\{
\begin{array}{ll}
\Psi_j'(s), &\mbox{ if there exists }\Psi_j'(s),
\\
0, & \mbox{ otherwise, }
\end{array}
\right. $ \;\; here $j=1,2$.\hfill \ensuremath{\triangleright}
\end{definition}
\begin{lemma}\label{zv29} 
If $\varkappa>0$, then {\it on some probability space}
there exists two r.v. $\xi_1'$ and $\xi_2'$ such, that:
\begin{equation}\label{zv14}
\mathsf{P} \{\xi_j'\leqslant s\}=\mathsf{P} \{\xi_j\leqslant s\}=\Psi_j;
\end{equation}
\begin{equation}\label{zv15}
\mathsf{P} \{\xi_1'=\xi_2'\}\geqslant \varkappa.
\end{equation}
\hfill \ensuremath{\triangleright}
\end{lemma}

\textbf{Proof.}
Denote
$$\begin{array}{ll}
\psi(u)\stackrel{\text{\rm def}}{=\!=} \min(\psi_1(u),\psi_2(u));&
\Psi(s)\stackrel{\text{\rm def}}{=\!=} \displaystyle
\int\limits_{-\infty}^s \psi(u)\,\mathrm{d}\, u;
\\
\hat{\Psi}(s)\stackrel{\text{\rm def}}{=\!=} \displaystyle
\frac{\Psi(s)}{\varkappa}; & \hat{\Psi}_j(s)\stackrel{\text{\rm
def}}{=\!=} \displaystyle \frac{\Psi_j(s)-\Psi(s)}{1-\varkappa}, \;
j=1,2
\end{array}
$$
(if $\varkappa=1$, then $\hat{\Psi}_j(s)\equiv 0$).

Let $\mathscr{U} _1$, $\mathscr{U} _2$ and $\mathscr{U} _3$ be independent r.v. with the continuous uniform distribution on
$[0;1]$.
Put $ \displaystyle \xi_j'\stackrel{\text{\rm
def}}{=\!=} \hat{\Psi}^{-1}(\mathscr{U} _j)\times \mathbf{1}
(\mathscr{U} _3\leqslant\varkappa) + \hat{\Psi}^{-1}_j(\mathscr{U}
_j)\times \mathbf{1} (\mathscr{U} _3>\varkappa) $.

Recall that the inverse function of monotonous function $f(s)$ is $
f^{-1}(u)\stackrel{\text{\rm def}}{=\!=} \inf\{s:\;f(s)\geqslant
u\}. $

R.v. $\xi_j'$ satisfy conditions (\ref{zv14})--(\ref{zv15}).

Lemma is proved.
\hfill \ensuremath{\blacksquare}

\begin{remark}\label{zv27} 
It at some time $\vartheta$, $n_\vartheta=0$ (i.e. the queueing system is idle), then the residual time of the idle period of the process $X_t$ (described in the Section \ref{zv16}) has the density of distribution $ \displaystyle
\varphi(s)=\lambda(X_{\vartheta+s})\exp\left(\displaystyle
\int\limits_0^s \lambda(X_{\vartheta+u})\,\mathrm{d}\,
u\right)\in(\lambda_0e^{-\Lambda s}; \Lambda_0e^{-\lambda s}). $
Thus, if two processes $X_t$ and $X_t'$ with different initial states $X_0$ and $X_0'$, at some time $\vartheta$ are in the set $\mathscr{S} _0$ (the queueing system is idle), then the common part $\varkappa_\vartheta$ of the distributions of their residual times of the staying in the set $\mathscr{S} _0$ satisfies the inequality
\begin{equation}\label{zv28} 
\widehat\varkappa_\vartheta\geqslant \int\limits_0^\infty \lambda_0e^{-\Lambda s}\,\mathrm{d}\, s=\frac{\lambda_0}{\Lambda}\stackrel{\text{\rm def}}{=\!=} \widehat\varkappa.
\end{equation}
\hfill \ensuremath{\triangleright}
\end{remark}

\section{Construction of successful coupling for two pro\-ces\-ses $X_t$, $X_t'$}
Consider two independent processes: the process $X_t$ with initial state $X_0=(0,0)$ and the process $X_t'$ with initial state $
\displaystyle X_0'=\linebreak\left(n_0',x'_0{}^{(0)};
x'{}^{(1)}_0,x_0'{}^{(2)},\ldots,x_0'{}^{(n_0')}\right). $

Their successful coupling $\mathscr{Z} _t=(Y_t,Y_t')$ will be constructed on the probability space $(\Omega,\mathscr{F} , \mathsf{P} )=\prod\limits_{i=0}^\infty (\Omega_i, \mathscr{F} _i , \mathsf{P} _i)$, where $\Omega_i=[0;1]$, $\mathscr{F} _i = \mathscr{B} ([0;1])$, and $\mathsf{P} _i $ is the Lebesgue measure on $[0;1]$.

The procedure of the construction of the successful coupling will be described {\it by steps}. This steps operate at the time moments when the first component of one of the processes 

$$Y_t=\left(m_t,y_t^{(0)};
y^{(1)}_t,y_t^{(2)},\ldots,y_t^{(m_t)}\right),\qquad
Y_t'=\left(m_t',y_t'{}^{(0)};
y_t'{}^{(1)},y_t'{}^{(2)},\ldots,y_t'{}^{(m_t')}\right)$$ 
changes.
Denote these time moments by $t_k$: $ \displaystyle t_0\stackrel{\text{\rm
def}}{=\!=} 0; \;\; t_k\stackrel{\text{\rm def}}{=\!=}
\inf\{t>t_{k-1}: m_{t}+m'_{t}\neq m_{t_{k-1}}+m'_{t_{k-1}} \}. $

In all steps, the residual \emph{virtual} times of the service of all customers being in the queueing systems (in hypothesis: during this service, the first component $n_t$ and $n_t'$ of processes not change), and the residual \emph{virtual} times until the arrival of the next customers.
For this aim, in any step, some quantity $\nu_k$ of exemplares $(\Omega_i,
\mathscr{F} _i , \mathsf{P} _i)$ is used, and in any step, new exemplares $(\Omega_i, \mathscr{F} _i , \mathsf{P} _i)$ are used.

{\bf Algorithm of construction}

\textbf{A. }If $m_{t_k}+m'_{t_k}>0$, for convenience put $Y_{\vartheta_k}=(\imath,\alpha_0;\alpha_1, \alpha_2, \ldots \alpha_\imath)$, $Y_{\vartheta_k}'=(\jmath,\beta_0;\beta_1, \beta_2, \ldots \beta_\jmath)$.

In $k$-th step, the probability spaces $(\Omega_i,
\mathscr{F} _i , \mathsf{P} _i)$, $i=\ell_k+1, \ell_k+2,\ldots
\ell_{k+1}$ are used; $\nu_k\stackrel{\text{\rm def}}{=\!=}
\ell_{k+1}-\ell_k=n_{\vartheta_k}+ n_{\vartheta_k}'+2=
\imath+\jmath+2$.
So, before $k$-th step, $\ell_k=\displaystyle \sum\limits _{i=0}^{k-1}\nu_i$ exemplars of the space $(\Omega_i, \mathscr{F} _i , \mathsf{P} _i)$ was used.

At the time $t_k$ \emph{virtual }time until arrival of the next customer $\xi$ (given in this time no ends of service) for $Y_t$ has c.d.f.\\
$F_1(s)=1-\exp\left(\displaystyle
\int\limits_0^s\lambda(\imath,\alpha_0+u;\alpha_1+u, \alpha_2+u,
\ldots \alpha_\imath+u)\,\mathrm{d}\, u\right)$.

R.v. $\xi\stackrel{\text{\rm def}}{=\!=}
F_1^{-1}(\omega)$, $\omega\in \Omega_{\ell_k+1}$ is defined on the space $(\Omega_{\ell_k+1}, \mathscr{F} _{\ell_k+1} ,
\mathsf{P} _{\ell_k+1})$.

The  \emph{virtual } residual time $\eta_j$ of the service of
$j$-th ($j=1,2,\ldots,\imath$) customer (given in this residual time no ends of other services end no new customers) for $Y_t$ has c.d.f. \\
$G_{1,j}(s)=1-\exp\left(\displaystyle \int\limits_0^s
h_j(\ell,\alpha_0+u;\alpha_1+u, \alpha_2+u, \ldots
\alpha_\imath+u)\,\mathrm{d}\, u\right)$ at the time $t_k$.

R.v. $\eta_j\stackrel{\text{\rm def}}{=\!=} G_{1,j}^{-1}(\omega)$,
$\omega\in \Omega_{\ell_k+j+1}$ is defined on the space $(\Omega_{\ell_k+j+1}, \mathscr{F}
_{\ell_k+j+1} , \linebreak\mathsf{P} _{\ell_k+j+1})$.
.

The construction of $\xi'$ and $\eta_i'$ is analogous
($i=1,2,\ldots,\jmath)$:
\\
$F_2(s)=1-\exp\left(\displaystyle
\int\limits_0^s\lambda(\jmath,\beta_0+u;\beta_1+u, \beta_2+u, \ldots
\beta_\jmath+u)\,\mathrm{d}\, u\right)$, $\xi'\stackrel{\text{\rm
def}}{=\!=} F_2^{-1}(\omega)$, $\omega\in \Omega_{\ell_k+\imath
+1}$;
\\
$G_{2,j}(s)=1-\exp\left(\displaystyle \int\limits_0^s
h_j(\jmath,\beta_0+u;\beta_1+u, \beta_2+u, \ldots
\beta_\jmath+u)\,\mathrm{d}\, u\right)$, $\eta_j'\stackrel{\text{\rm
def}}{=\!=} G_{2,j}^{-1}(\omega)$, $\omega\in
\Omega_{\ell_k+\imath+j+2}$.

In $k$-th step, $\nu_k=\imath+\jmath+2$ exemplares $(\Omega_i, \mathscr{F} _i , \mathsf{P} _i)$ was used; $\ell_{k+1}=\ell_k+\nu_k$.

The time $t_{k+1}=t_k+\min\{\xi,\xi',\eta_j,\eta'_i \;
(j=1,2,\ldots \imath,\; i=1,2,\ldots \jmath)\}$ is the time of the change of the first component of the process $Y_t$or of the process $Y_t'$.

If again $m_{t_k}+m'_{t_k}>0$, then the procedure \textbf{A} is repeated.

\textbf{B.} If $m_{t_k}+m'_{t_k}=0$, then $Y_{t_k}=(0,\alpha_0)$ and
$Y'_{t_k}=(0,\beta_0)$. C.d.f. of the times $\xi$ and $\xi'$ till the arrival of the new customer for the process $Y_t$ or $Y_t'$ correspondingly, is
\begin{eqnarray*}
\mathsf{P} \{\xi\leqslant s\}=F_1(s)=1-\exp\left(\displaystyle \int\limits_0^s\lambda(0,\alpha_0+u)\,\mathrm{d}\, u\right);
\\
\mathsf{P} \{\xi'\leqslant s\}=F_2(s)=1-\exp\left(\displaystyle \int\limits_0^s\lambda(0,\beta_0+u)\,\mathrm{d}\, u\right).
\end{eqnarray*}
By Remark \ref{zv27} and (\ref{zv28}), the common part of these distributions is not less than $\widehat\varkappa=\displaystyle
\frac{\lambda_0}{\Lambda}$.
Therefore, on the space $\displaystyle
\prod\limits_{i=\ell_k+1} ^{\ell_k+3} (\Omega_i, \mathscr{F} _i ,
\mathsf{P} _i)$ (here $\nu_k=3$) By Lemma \ref{zv29}
there exists $\xi$ and $\xi'$ such that $\mathsf{P}
\{\xi=\xi'\}\geqslant \widehat\varkappa$.

At the time moment $t_{k+1}\stackrel{\text{\rm def}}{=\!=}
t_k+\min\{\xi,\xi'\}$ the first component of at least one of the processes
$Y_t$, $Y_t'$ increases, and the first components of both processes $Y_t$, $Y_t'$ will change simultaneously -- with probability greater than $\widehat\varkappa$,
i.e. the event $
Y_{t_{k+1}}=Y_{t_{k+1}}'=(1,0;0)$ will happen.

\textbf{C. } If the event $Y_{t_{k}}=Y_{t_{k}}'=(1,0;0)$ happened, then for $t>t_k$ put the processes $Y_t$ and $Y_t'$ equal. Namely, it can construct only one process $Y_t$ by procedure given in \textbf{A}. And $Y_t'\equiv Y_t$ for $t>t_{k+1}$.

The processes $Y_t$ and $Y_t'$ constructed by the procedures \textbf{A-C} satisfy the conditions (i)--(ii) from Definition \ref{zv30}.

Let us check the condition (iii) from Definition \ref{zv30}.

Consider the process $Y_t'\stackrel{\mathscr{D}}{=} X_t'$.
If $n_0'>0$, then by comparison with the ``classic'' queueing system (see Sections
\ref{zv20} and \ref{zv31}), it can estimate the moments of the residual time $\zeta_1$ of the first busy period of $Y'_t$; it is easy to see that for $k\in[0;K)$, $\mathbb E (\zeta_1)^k<\infty$.

Denote $\vartheta_0\stackrel{\text{\rm def}}{=\!=} \zeta_0$,
$\vartheta_k\stackrel{\text{\rm def}}{=\!=} \inf\{t>\vartheta_k:\;
Y'_t\in\mathscr{S} _0, Y'_{t-0}\notin\mathscr{S} _0\}$;
$\vartheta_k$ is the time of end of $k$-th busy period of
the process $Y_t'$.
If $Y_0'\in\mathscr{S} $, put
$\vartheta_0\stackrel{\text{\rm def}}{=\!=} 0$.

As distribution the process $Y_t$ is equal to distribution of the
the process $X_t$, the inequality $\mathsf{P}
\{X_{\vartheta_k}\in\mathscr{S} _0\}\geqslant e^{-\rho}$ is true -- see
(\ref{zv32}).

In addition, if both processes $Y_t$ and $Y_t'$ hit to the set $\mathscr{S} $, then at the exit from this set, they coincide with probability greater then $\widehat\varkappa$.
Hence, again denote $\tau\stackrel{\text{\rm def}}{=\!=} \inf\{t>0:\;Y_t=Y_t'\}$, and $\eta_k$ is the length of an idle period of $Y_t'$ after $\vartheta_k$, and we have:
$ \displaystyle
\mathsf{P} \{\tau>\vartheta_k+\eta_k\}\leqslant (1-\varpi)^k.
$
Note, that $\tau=\tau(Y_0')=\tau(X_0')$ (as $X_0=Y_0=(0,0)$ is fixed).

Evidently, $\vartheta_k+\eta_k=\vartheta_0+\eta_0+\displaystyle \sum\limits _{i=1}^k ( \zeta_i+\eta_i)$, where $\zeta_i$ is $i$-th busy period, and $\eta_i$ is the next free period.
Denote $\mathscr{E} _k\stackrel{\text{\rm def}}{=\!=} \{Y_{\vartheta_k}\neq Y_{\vartheta_k}'\;\&\;Y_{\vartheta_k+\eta_k}= Y_{\vartheta_k+\eta_k}' \}$ -- it is the event ``coincidence of the processes $Y_t$ and $Y_t'$ is immediately after $k$-th hit of $Y_t'$ t the set $\mathscr{S} _0$''.

\begin{proposition} \label{zv34} 
If $r\in[1; K)$ then
$$\displaystyle
\mathbb{E} \,\tau^r\leqslant C_1(\varpi,r) \mathbb{E} \, \vartheta_0^r+ C_1(\varpi,r)\mathbb{E} \,\zeta_1^r+ C_3 (\varpi,r,\lambda_0,\Lambda),
$$
where the constatns $C_1$, $C_2$, $C_3$ can be calculated by the formulae below. \hfill \ensuremath{\triangleright}
\end{proposition}

\textbf{Proof.} Firstly, estimate the bounds for $\mathbb{E}
\,(\eta_k)^r|_{\mathscr{E} _k}$. If the event $\mathscr{E} _k$ happened, then 
r.v. $\eta_k$ (see Lemma \ref{zv29}) has a density of distribution
\begin{multline*}
\varphi(s)= \min\left\{\lambda (Y_{\vartheta_k}+s)\exp \left(-\int\limits_0^ s\lambda(Y_{ \vartheta_k}+u)\,\mathrm{d}\, u\right),\right.
\\
\left.\lambda (Y'_{ \vartheta_k}+s)\exp \left(-\int\limits_0^s \lambda(Y'_{ \vartheta_k}+u)\right) \right\}
\leqslant\Lambda e^{-\lambda_0 s}.
\end{multline*}
Hence, $\mathbb{E} \,(\eta_k)^r|_{\mathscr{E} _k}\leqslant\displaystyle \frac{r! \Lambda}{\lambda_0^{r+1}}$.
Also $\mathbb{E} \,(\eta_k)^r|_{\overline{\mathscr{E} }_k}\leqslant\displaystyle \frac{r! \Lambda}{\lambda_0^{r+1}}$; $\zeta_k$ not depends on $\mathscr{E} _j$ if $k\neq j$.

Now for estimation of $\mathbb{E} \,\tau^r$ we use the Jensen's inequality in the form
(\ref{zv33})):
\begin{multline*}
\mathbb{E} \,\tau^r\leqslant \displaystyle \sum\limits _{k=1}^\infty \mathbb{E} \left(\vartheta_0+\eta_0+\displaystyle \sum\limits _{i=1}^k ( \zeta_i+\eta_i)\right)^r\mathbf{1} (\mathscr{E} _k)\leqslant
\\
\leqslant\displaystyle \sum\limits _{k=1}^\infty (1-\varpi)^{k-1}\varpi (2k+2)^{r-1}\left( \left(\mathbb{E} \,\vartheta_0^r+(k+1)\displaystyle \frac{r! \Lambda}{\lambda_0^{r+1}}+k\mathbb{E} \,\zeta_1^r \right)\right)<
\\
<2^{r-1}\left(\mathbb{E} \,\vartheta_0^r\displaystyle \sum\limits _{k=1}^\infty (1-\varpi)^{k-1}(k+1)^{r-1}+ \displaystyle \frac{r! \Lambda}{\lambda_0^{r+1}} \displaystyle \sum\limits _{k=1}^\infty (1-\varpi)^{k-1}(k+1)^r+\right.
\\
+\left.\mathbb{E} \,\zeta_1^r\displaystyle \sum\limits _{k=1}^\infty (1-\varpi)^{k-1}(k+1)^{r-1}k\right).
\end{multline*}
Proposition \ref{zv34} is proved.
\hfill \ensuremath{\blacksquare}

So now, the constructed pair $\mathscr{Z} _t=(Y_t,Y_t')$ is a successful coupling of the processes $X_t$ and $X_t'$, and there is the explicit bounds for the coupling time:
\begin{equation}\label{zv39} 
\mathbb{E} \,(\tau(X_0'))^r\leqslant \mathcal{C}(r,\lambda_0, \Lambda,\mathbb{E} \,\vartheta_0^r).
\end{equation}

Therefore,
\begin{equation}\label{zv40} 
|\mathsf{P} \{X_t\in S\}-\mathsf{P} \{X_t'\in S\}|\leqslant \displaystyle \frac{\mathcal{C}(r, \lambda_0,\Lambda,\mathbb{E} \,\vartheta_0^r)}{t^r},
\end{equation}
but here the next stage -- integration with respect to the stationary measure (\ref{zv35}) -- is impossible, because the stationary distribution of the queueing system $M|G|\infty$ (the process $\stackrel{\circ}{X}_t$) is \emph{unknown}, and for the studied process ${X}_t$ also.

\section{Process $X_t$ as a regenerative process}

\begin{definition}\label{rege} 
The stochastic process $\{W_t,\,t\geqslant0\}$, defined on the probability space $(\Omega, \mathscr F, \mathbf{P})$, measurable with respect to the filtration $\{\mathscr{F}_t,\,t\geqslant 0\}$, with the state space
$(\mathscr{W},\mathscr{B}(\mathscr{W}))$, is called 
\textit{regenerative process}, if there exists the sequence of the Markov times (stopping times) (with respect to filtration $\mathscr{F}_t$) $\{\theta_i\}_{i\in \mathbb N}$, such that
\begin{enumerate}
\item $W_{\theta_i}=W_{\theta_j}$ for all $i,j\in\mathbb N$;
\item Random elements $\Xi_i \stackrel{{\rm def}}{=\!\!\!=} \{W_t, \, t\in[\theta_i,\theta_{i+1}]\}$ ($i\in\mathbb N$) are i.i.d.\hfill \ensuremath{\triangleright}
\end{enumerate}
\end{definition}
The intervals $[\theta_i,\theta_{i+1}]$ are called regeneration periods.
Evidently, the length of the regeneration periods $\chi_i\stackrel{\text{\rm
def}}{=\!=} \theta_{i+1}-\theta_i$ are i.i.d.. 
Denote $\Phi(s)\stackrel{\text{\rm
def}}{=\!=} \mathsf{P} \{\chi_i\leqslant s\}$ -- c.d.f. of r.v. $\chi_i$; it is unknown.

If $\mathbb{E} \,\chi_i< \infty$, then the distribution of the regenerative process $W_t$ weakly converges to the stationary (invariant) distribution -- it is well-known fact.

Note, for the stochastic process $X_t=\left(n_t, x^{(0)}_t; x^{(1)}_t,x^{(2)}_t,\ldots,x^{(n_t)}_t\right)$ described in Section \ref{zv16}, the regeneration periods are the sum of busy period and next idle period: $\chi_k=\zeta_k+\eta_k$.

The values of $\mathbb E (\zeta_k)^r$ and $\mathbb E (\eta_k)^r$ can be estimated by (\ref{zzzv}), (\ref{zv36}) and (\ref{zv25}).

Note, the regeneration period is greater then its part -- idle period having the distribution density $\varphi(s)\geqslant \lambda_0 e^{-\Lambda s}$ (Remark \ref{zv27}).

So, the next bounds are true:
\begin{equation}\label{zv43} 
\begin{array}{l}
\mathbb{E} \,\chi_k^r\leqslant 2^{r-1}(\mathbb{E} \,\zeta_1^r+\mathbb{E} \,\eta_1^r)=\mathcal{M}_r,\qquad
\mathbb{E} \,\chi_i>\displaystyle \frac{\lambda_0}{\Lambda^2}.
\end{array}
\end{equation}

The Markov times $\{\theta_i\}_{i\in \mathbb N}$ form the embedded renewal process \linebreak $N_t=\displaystyle
\sum\limits _{i=1}^\infty \mathbf{1} (\theta_i<t)$.

Consider the \emph{backward renewal time} of the process $N_t$. 
This is Markov processs $B_t\stackrel{\text{\rm def}}{=\!=} t-\theta_{N_t}$ (with the state space $\mathbb{R}_+$).
It is well-known (\cite{GBS,GK}), that the distribution of the backward renewal time $\mathfrak{P}_t$ weakly converges to the stationary invariant distribution $\mathfrak{P}$ with c.d.f. \begin{equation}\label{zv37} 
\widetilde{\Phi}(s)=\displaystyle \frac{\displaystyle \int\limits_0^s(1-\Phi(u))\,\mathrm{d}\, u}{\mathbb{E} \,\chi_i}.
\end{equation}
Also, this distribution (\ref{zv37}) is a stationary distribution for the \emph{forward renewal process}
$R_t\stackrel{\text{\rm def}}{=\!=}
\theta_{N_t+1}-t$ (here the renewal processes $N_t$).

As the process $X_t$ is regenerative, its distribution fully defined by the time elapsed from the last regeneration point (beginning of the regeneration period).
So, $\mathscr{P} _t=\mathbb{F} (B_t)$.

Now, consider the process $X_t$ with initial state $(0,0)$ and the process $X_t'$ with arbitrary initial state $X_0$.
Both these processes have the embedded renewal processes: $N_t$ and $N_t'$ correspondingly. Moreover, for the renewal processes $N_t$ and $N_t'$ the backward renewal processes $R_t$ and $R_t'$ correspond.

The construction of the successful coupling $\mathscr{Z} _t=(Y_t,Y'_t)$ for the pair of the processes $X_t$ and $X_t'$, gives simultaneously the successful coupling $\mathscr{W} _t=(D_t,D_t')$ for the pair of the processes $B_t$ and $B_t'$.
Denote $\mathfrak{P}_t$ and $\mathfrak{P}_t'$ -- the distributions
of the processes $B_t$ and $B_t'$ correspondingly ($\mathsf{P} \{B_t\in
S\}=\mathfrak{P}_t(S)$ and $\mathsf{P} \{B_t'\in
S\}=\mathfrak{P}_t'(S)$).

From (\ref{zv39}) and (\ref{zv40}) we have
\begin{equation}\label{zv41} 
|\mathsf{P} \{B_t\in S\}-\mathsf{P} \{B_t'\in S\}|\leqslant \displaystyle \frac{\mathcal{C}(\varpi,r, \lambda_0,\Lambda,\mathbb{E} \,\vartheta_0^r)}{t^r},
\end{equation}
where $\vartheta_0=R_0'$ is the backward renewal process of the renewal process $N_t'$ at the time moment $t=0$.
\begin{proposition}\label{zv42} 
\begin{multline*}
\|\mathfrak{P}_t-\mathfrak{P}\|_{TV}\leqslant
\\
\leqslant \frac{2\left( C_1(\varpi,r) \displaystyle \frac{\mathbb{E} \, \chi_i^{r+1}}{(r+1)\mathbb{E} \,\chi_i}+ C_1(\varpi,r)\mathbb{E} \,\zeta_1^r+ C_3 (\varpi,r,\lambda_0,\Lambda)\right)} {t^r} = \displaystyle \frac{\mathcal{C}_1(r,\lambda_0,\Lambda)}{t^r},
\end{multline*}
where $r\in[1,K-1)$. \hfill \ensuremath{\triangleright}
\end{proposition}
\textbf{Proof.}
Consider the backward renewal process $B_t$ with the initial state $B_0=0$ and the backward renewal process $\widetilde{B}_t$ with the initial stationary distribution $\mathfrak{P}$; for all $t\geqslant 0$
\begin{equation}\label{zv44} 
\mathsf{P} \{\widetilde{B}_t\in S\}=\mathfrak{P}(S)=\int\limits_{S}\frac{1-\Phi(u)\,\mathrm{d}\, u}{\mathbb{E} \,\chi_i}.
\end{equation}

By substitution in (\ref{zv41}) $\vartheta_0=R_0'$ by the stationary forward renewal process $\widetilde{R}_0$ with the distribution (\ref{zv37}), the Proposition \ref{zv42} is proved.

Considering (\ref{zv43}), it can calculate the upper bound for the value $\mathcal{C}_1(r,\lambda_0,\Lambda)$.
\hfill \ensuremath{\blacksquare}

\begin{theorem} In the conditions (\ref{zv10}), the inequality
$\displaystyle
\|\mathscr{P} _t-\mathscr{P} \|_{TV}\leqslant \displaystyle \frac{\mathcal{C}_1(r,\lambda_0,\Lambda)}{t^r}
$
is true for all $r\in[1,K-1)$.
\hfill \ensuremath{\triangleright}
\end{theorem}
\textbf{Proof.} It was said above that $\mathscr{P}
_t=\mathbb{F}(\mathfrak{P}_t)$. 
Accordingly, $\mathscr{P}
=\mathbb{F}(\mathfrak{P})$. 
Recall that $\chi_i$ is the length of the regeneration period, which includes absolutely continuous r.v. So, the distribution of r.v. $\chi_i$ is absolutely continuous. 
Therefore, the distribution of the backward renewal process $B_t$ has a density of c.d.f
$\psi_t(s)$.
It is known (\ref{zv44}), the c.d.f. of the stationary distribution of the backward renewal process is known, its density is $\psi(s)=\displaystyle \frac{1-\Phi(u)}{\mathbb{E} \,\chi_i}$.

Denote $\sigma_t=\{s: \psi_t(s)> \psi(s)\}$.

Now,
\begin{multline*}
|\mathscr{P} _t(S)-\mathscr{P} (S)|=
\\
=\left|\int\limits_0^\infty \mathsf{P} \{X_t\in S|B_t=u\}\psi_t(u)\,\mathrm{d}\, u - \int\limits_0^\infty \mathsf{P} \{X_t\in S|B_t=u\}\psi(u)\,\mathrm{d}\, u\right|\leqslant
\\
\leqslant \max \left\{\int\limits_{\sigma_t}(\psi_t(u)-\psi (u))\,\mathrm{d}\, u, \int\limits_{\mathbb{R}_+\setminus \sigma_t}(\psi(u)-\psi_t (u))\,\mathrm{d}\, u \right\}=
\\
=\max\left\{|\mathfrak{P}_t(\sigma_t)- \mathfrak{P}(\sigma_t))|, |\mathfrak{P}_t(\mathbb{R}_+\setminus\sigma_t)- \mathfrak{P}(\mathbb{R}_+ \setminus\sigma_t)|\right\} \leqslant
\\
\leqslant \sup\limits _{\sigma\in \mathscr{B} (\mathbb{R})}| \mathfrak{P}_t(\sigma) -\mathfrak{P}(\sigma)|,
\end{multline*}
so $\|\mathscr{P} _t-\mathscr{P} \|_{TV}\leqslant \|\mathfrak{P}_t-\mathfrak{P}\|_{TV}$, and it is the end of the Proof.
\hfill \ensuremath{\blacksquare}

\begin{remark}\label{fin} 
The obtained bounds are not optimal for the following reasons.
\begin{enumerate}
\item The application of the coupling method to processes in continuous time usually gives a suboptimal estimates.
\item In the construction of the successful coupling only times of hitting of the process $X_t'$ to the idle state has been used.
In addition, it can consider the times of the hitting of the process $X_t$ to the idle state.
\item In the conditions (\ref{zv10}), it is impossible to give the accurate estimate for ($\mathcal{C}$ and $\mathcal{C}_1$).
If the distributions or intensities are given more accurate, these values can be estimated more exactly.
\end{enumerate}
However, the resulting estimates may be useful in practical applications, because the simulation is not always able to provide sufficient qualitative bounds.
\hfill \ensuremath{\triangleright}
\end{remark}

\paragraph{Acknowledgement.} The author thanks L. G. Afanasyeva, A. Yu. Veretennikov and A. D. Manita for important and useful comments in the discussion of this work. Also, the author thanks E. Yu. Kalimulina for the great help in the preparation of the text.
The work is supported by RFBR (project No 17-01-00633 A).

\end{document}


\section*{Об одной обобщённой системе Эрланга-Се\-вас\-тья\-но\-ва и её скорости сходимости.}
\begin{flushright}
\textbf{Г. А. ЗВЕРКИНА}\\

{\it  Российский университет транспорта (МИИТ),\\
Институт проблем управления\\
им. В.А.Трапезникова РАН\\
zverkina@gmail.com\\}
\end{flushright}
УДК 519.21\\
{\small \textbf{Ключевые слова:} распределение состояния системы массового обслуживания, скорость сходимости, регенерирующие процессы, метод склеивания.
\begin{center}
\textbf{Аннотация}
\end{center}
Получена оценка сверху для скорости сходимости распределения состояния бесконечнолинейной системы массового обслуживания в случае, когда интенсивность входящего потока и интенсивность обслуживания зависят от состояния системы.}

{\small \begin{center}
\textbf{Abstract}
\end{center}
{\it G. A. Zverkina. About some extended Erlang-Sevastyanov queueing system and its con\-ver\-gen\-ce rate. Fundamentalnaya i
prikladnaya matematika, vol. **, no. *, pp. ** - **}

The upper bound for the convergence rate of the distribution  of a queuing system state with infinitely many servers is obtained for the case when the intensities of the incoming  and service flows depend on the state of the system. }

\section{Введение. }
\subsection{Мотивация}

Рассматривается бесконечнолинейная система массового обслуживания (СМО), в которую в моменты времени $t_1$, $t_2$,\ldots, $t_n$ поступают заявки, которые сразу же начинают обслуживаться (полагаем $t_0\stackrel{\text{\rm def}}{=\!=}  0$). Время обслуживания $i$-й заявки -- случайная величина (сл.в.) $\xi_i$ с функцией распределения (ф.р.) $G_i(s)$, а промежутки между поступлениями требований $\eta_i
\stackrel{\text{\rm def}}{=\!=}  t_i-t_{i-1}$ -- сл.в. с ф.р. $F_i(s)$.
Для простоты предположим, что распределения сл.в.  $\xi_i$ и $\eta_i $ абсолютно непрерывны; обозначим $f_i(s)\stackrel{\text{\rm def}}{=\!=}  F'_i(s)$, $g_i(s)\stackrel{\text{\rm def}}{=\!=}  G'_i(s)$.

Мы будем описывать распределения $F_i(s)$ и $G_i(s)$ в терминах {\it интенсивностей}.
\begin{definition}\label{zv1}  
{\it Интенсивностью} в момент времен $s>0$ окончания временн\'ого периода $\mathscr{T}$, начавшегося в момент времени $t=0$ и не завершившегося до момента времени $s>0$,  называется
$
\displaystyle\phi(s)\stackrel{\text{\rm def}}{=\!=}  \lim\limits _{\Delta \downarrow 0}\frac{\mathsf{P} \{\mathscr{T} \in (s,s+\Delta]| \mathscr{T}>s\}}{\Delta}.
$
Если $\mathscr{T}$ имеет абсолютно непрерывную ф.р. $\Phi(x)$, то

\hspace{2.5cm}$
\displaystyle\phi(s)=  \lim\limits _{\Delta \downarrow 0} \frac{\Phi(s+\Delta)-\Phi(s)}{\Delta(1-\Phi(s))}= \frac{\Phi'(s+0)}{1-\Phi(s)}.
$\hfill \ensuremath{\triangleright}
\end{definition}

\begin{remark}
Из Определения \ref{zv1} следует, что
$
\mathsf{P} \{\mathscr{T}\in(s,s+\Delta]| \mathscr{T}>s\}=\phi(s)\Delta+ o(\Delta).
$\hfill \ensuremath{\triangleright}
\end{remark}
\begin{remark}\label{zv9}  
Функция распределения и её плотность легко  восстанавливаются по интенсивности:
\begin{equation}\label{zv2}  
    \Phi(s)=1-\exp\left(-\int\limits_0^s \phi(u)\,\mathrm{d}\, u\right);
\end{equation}
\begin{equation}\label{zv3}  
\Phi'(s)= \phi(s)\exp\left(-\int\limits_0^s \phi(u)\,\mathrm{d}\, u\right).
\end{equation}
Очевидно, постоянная интенсивность соответствует экспоненциальному распределению.
Интенсивность можно также рассматривать и в случае, когда ф.р. $\Phi(s)$ имеет разрывы. Но в этом случае формулы (\ref{zv2}) и  (\ref{zv3}) неприменимы.
\hfill \ensuremath{\triangleright}
\end{remark}
Итак, мы описываем поведение бесконечнолинейной СМО на языке интенсивностей, а именно.
\\
{\bf 1.} Интенсивность поступления очередной заявки в момент времени $t$ есть $\lambda(t)$, т.е.
$$
\mathsf{P} \left\{
\begin{array}{l}
\mbox{Очередная заявка поступит}
\\
\mbox{в интервале $(t, t+\Delta]$}
\end{array}
\right \}=\lambda(t)\Delta+ o(\Delta).
$$
\\
{\bf 2.} Интенсивность обслуживания $i$-й по порядку поступления заявки из имеющихся $n$ заявок в системе  в момент времени $t$ равна $h_i(t)$, т.е.
\begin{multline*}
\mathsf{P} \left\{
\begin{array}{l}
\mbox{$i$-е требование, обслуживающееся }
\\
\mbox{в момент времени $t$, }
\\
\mbox{закончит обслуживание }
\\
\mbox{на промежутке $(t,t+\Delta]$}
\end{array}
\right\}=
\\
= \mathsf{P} \{\xi_i\in(s,s+\Delta]|\xi_i>s\}=h_i(t)\Delta+ o(\Delta).
\end{multline*}

{\it Здесь и везде далее под $i$-м требованием, находящимся в СМО, мы понимаем $i$-е по порядку поступления требование из находящихся в СМО.}

В случае, когда $\lambda(t)\equiv \lambda$ и $h_i(t)\equiv \mu$, мы  имеем дело с СМО $M|M|\infty$, впервые изученной в \cite{Zv8}; было показано, что при условии $\mu<\lambda$ распределение количества заявок в системе $n_t$ сходится к стационарному распределению:

\hspace{2cm}$
\displaystyle\lim\limits _{t\to\infty}\mathsf{P} \{n_t=k\}=\mathscr{P} _k=\mathscr{P} _0\frac{\lambda^k}{\mu^k}, \mbox{ где } \mathscr{P} _0=1-\frac{\lambda}{\mu}.
$

Позднее этот результат был обобщен на случай, когда интенсивность входящего потока зависит только от числа требований в системе, т.е. $\lambda(t)\equiv \lambda_k$ при $n_t=k$ (см., например, \cite[Гл.4, \S 4-5]{zv4}):
\begin{equation}\label{zv4} 
\lim\limits _{t\to\infty}\mathsf{P} \{n_t=k\}=\mathscr{P} _k=\displaystyle \mathscr{P} _0\frac{\displaystyle \prod\limits_ {i=0} ^{k-1}\lambda_i}{\mu^k}, \mbox{ где } \mathscr{P} _0=\left(\sum\limits _{k=0}^\infty \frac{\displaystyle \prod\limits_ {i=0} ^{k-1}\lambda_i} {\mu^k}\right)^{-1}.
\end{equation}
Более того, если интенсивность обслуживания всех находящихся в СМО заявок зависит только от количества заявок, т.е. $h_i(t) \equiv \mu_k>0$ при  $n_t=k$, формула \ref{zv4} применима, но с заменой $\mu^k$ на $\prod\limits_ {i=1} ^{k}\mu_i$.
Очевидно, если $\lambda_k=0$ при некотором $k$, то (\ref{zv4}) применима для системы с конечным числом серверов без очереди.

В описанном случае, когда
\begin{equation}\label{zv5}  
\begin{array}{l}
\lambda(t)\not\equiv \lambda,\qquad \mu(t)\not\equiv \mu,
\\
\lambda(t)=\lambda_{n_t}, \qquad h_i(t)=\mu_{n_t},
\end{array}
\end{equation}
где  $n_t$ -- количество требований в системе в момент времени $t$, времена обслуживания и входящий поток зависимы между собой.

В \cite{Zv11} исследовалась СМО с пуассоновским входящим потоком и произвольным абсолютно непрерывным распределением времени обслуживания (СМО $M|G|\infty$), одинаковым для всех требований; предполагалось, что времена обслуживания и входящий поток независимы.
В терминах интенсивностей это означает, что если в момент времени $t$ в системе находится $n_t$ требований, и $i$-е требование  находилось на обслуживании в течении времени $x^{(i)}_t$ ($i=1,2,\ldots,n$), то интенсивность обслуживания $i$-го требования в момент времени $t$ равна $h_i(t)=h\left(x^{(i)}_t\right)$.

Обозначим $\mu\stackrel{\text{\rm def}}{=\!=}  \left(\displaystyle \int\limits_0^\infty s\,\mathrm{d}\, G(s)\right)^{-1}$, где  $G(s)=1-\exp\left(-\displaystyle \int\limits_0^s h(u)\,\mathrm{d}\, u\right)$ -- ф.р. времени обслуживания.
В \cite{Zv11} было высказано утверждение, что для $k>0$
\begin{multline*}
\lim\limits _{t\to\infty} \mathsf{P} \{n_t=k; x_t^{(1)}<\alpha_1, x_t^{(2)}<\alpha_2,\ldots x_t^{(k)}<\alpha_k  \}=
\\
=\mathscr{P} (k;\alpha_1, \alpha_2, \ldots, \alpha_k)=\mathscr{P} _0\prod\limits_{i=1}^k \lambda_i G(1-\alpha_i),
\end{multline*}
где нормирующий множитель $\mathscr{P} _0$ вычисляется так же, как в (\ref{zv4}).

Этот факт для СМО $M|G|m|0$ в более общем виде, без предположения об абсолютной непрерывности ф.р. времени обслуживания $G(s)$,  был доказан в \cite{seva1} (см. также \cite{seva2}).

То есть в \cite{Zv11} было предложено рассматривать СМО, в которой
\begin{equation}\label{zv6} 
\lambda(t)\equiv \lambda, \qquad h_i(t)=h_i\left(x_t^{(i)}\right),
\end{equation}
где $x_t^{(i)}$ -- прошедшее к моменту времени $t$ время обслуживания $i$-го по порядку поступления требования из $n_t$ требований, находящихся в системе.

В дальнейшем результаты, основанные на \cite{Zv8}, были получены многими исследователями (см., например, \cite{zv2,Zv9,Zv12,bor}).

Естественным обобщением ситуаций (\ref{zv5}) и  (\ref{zv6}) является СМО, в которой интенсивности входящего потока и обслуживания зависят не только от $n_t$ и прошедших времён обслуживания находящихся на обслуживании требований, но и от того, когда произошло последнее изменение $n_t$ -- количества требований в системе:

\hspace{0.3cm}$
\lambda(t)=\lambda\left(x_t^{(0)}; x_t^{(1)}, x_t^{(2)},\ldots,x_t^{(n_t)}\right), \;\; h_i(t)=h_i\left(x_t^{(0)}; x_t^{(1)}, x_t^{(2)},\ldots,x_t^{(n_t)}\right),
$
\\
где $x_t^{(0)}$ -- время, прошедшее с последнего изменения количества требований в системе $n_t$, а  $x_t^{(i)}$ -- прошедшее к моменту времени $t$ время обслуживания $i$-го по порядку поступления требования из $n_t$ требований, находящихся в системе.
Знание стационарного распределения состояния СМО чрезвычайно важно для приложений ТМО, на основании этих распределений рассчитываются технические характеристики сложных систем.

Такого рода СМО были рассмотрены, например, в \cite{ver1,ver2,zv7}.

Например, \cite{zv7} был доказан следующий факт:
для СМО с одним обслуживающим прибором и бесконечной очередью, в случае, когда интенсивность обслуживания $h$ и интенсивность входящего потока $\lambda$ зависят от количества заявок в системе $n_t$ и прошедшего времени обслуживания $x_t$, верна следующая
\begin{theorem}\label{zv7} 
Обозначим $X_t\stackrel{\text{\rm def}}{=\!=} (n_t,x_t)$, где $n_t$ -- число находящихся в СМО требований в момент времени $t$, $x_t$ -- время, в течении которого находилось на обслуживании требование, которое в момент времени $t$ находится на обслуживании.
В случае $n_t=0$ интенсивность входящего потока постоянна и равна $\lambda_0>0$ и $X_t\stackrel{\text{\rm def}}{=\!=}(0,0)$.

Если $\lambda(n_t,x_t)$ и $h(n_t,x_t)$ измеримы по Борелю и ограничены, выполнены условия
$$
\begin{array}{l}
\sup\limits _{n,x:\;n>0}\lambda(n,x)= \Lambda<\infty;\qquad \inf\limits _{n>0, x\geqslant 0} h(n,x)\geqslant\displaystyle  \frac{C_0}{1+x};
\\ \\
C_0>4(1+2\Lambda), \mbox{ и для некоторого } k>1 \quad C_0>2^{k+1}(1+\Lambda 2^k),
\end{array}
$$
то распределение процесса $X_t=(n_t,x_t)$ сходится к стационарному распределению $\mathscr{P} $ на пространстве состояний $\mathscr{X} = \{(0,0)\} \bigcup \{(n,x): \; n \in \mathbb N, \;
x\ge 0\}$, и  найдутся такое $m>k$ и $C>0$, что при всех $t\geqslant 0$
$$
\sup\limits _{S\in \mathscr{B} (\mathscr{X} )}|\mathsf{P} \{X_t\in S|X_0=(n,x)\}-\mathscr{P} (S)|\leqslant C\,\frac{(1+n+x)^m}{(1+x)^k},
$$
где $\mathscr{P} $ -- стационарное распределение процесса $X_t$.
\hfill \ensuremath{\triangleright}
\end{theorem}
В доказательстве этой Теоремы \ref{zv7} приводится (неполный) алгоритм вычисления константы $C$.

Однако существует достаточно большое количество результатов о {\it порядке} скорости сходимости распределения состояния СМО к стационарному распределению (см., например, \cite{ zv1,Zv14, GK,bor} и др.). В этих работах доказано, что некоторые характеристики СМО $\chi_t$ сходятся к стационарному значению $\tilde\chi$ со скоростью $\phi(t)$, т.е. существует такая константа $K$, что при всех $t$ выполняется 
$
|\chi_t-\tilde\chi|<K \phi(t),
$
где $\phi(t)\searrow 0$.

Так,  например, в СМО, исследованной в \cite{zv7} (Теорема \ref{zv7}), для сходимости распределения процесса $X_t$ к стационарному распределению со скоростью не превосходящей $\displaystyle \frac{K}{(1+x)^k}$, $k\in[0,C_0-1)$, достаточно  $\Lambda>C_0$, но без дополнительных условий Теоремы \ref{zv7} оценить константу  $K$ ранее не удавалось.

Отметим, что для СМО $MI|GI|\infty$ хорошо известны (см., например, \cite{Zv9}) такие характеристики, как $\mathscr{P} _k(t)$ -- вероятность того, что в системе занято $k$ приборов в момент времени $t$,  известны средняя длина периода занятости и его дисперсия, а также преобразование Лапласа ф.р. периода занятости (если в начальный момент времени система свободна).

Соответственно, эти результаты позволяют легко оценить  скорость сходимости величин $\mathscr{P} _k(t)$ к их стационарным значениям $\mathscr{P} _k$ -- см., например, \cite{Zv3}.

\subsection{Известные факты и цель работы}
Напомним, что СМО $MI|GI|\infty$ состоит из бесконечного числа приборов, работающих независимо, и случайное время обслуживания попавшего на прибор требования (заявки) одинаково распределено для всех приборов и имеет ф.р. $G(s)$; входящий поток пуассоновский с постоянным параметром $\Lambda$.
Предполагается, что в начальный момент времени система свободна (нет требований на обслуживании).

Цель настоящей работы -- исследование поведения бесконечнолинейной \linebreakСМО, в которой входящий поток и времена обслуживания задаются интенсивностями, зависящими от {\it полного состояния системы}; при этом  все времена обслуживания находящихся в системе требований {\it некоторым образом} зависят между собой, а также они зависимы и от входящего потока требований, который может менять свою интенсивность в зависимости от \emph{полного состояния системы}.
При этом естественно предполагать, что, несмотря на такую зависимость, некоторые характеристики  времени обслуживания и входящего потока находятся в определённых границах (например, естественно требовать, чтобы интенсивность входящего потока была ограничена, а интенсивность обслуживания не обращалась в нуль).

В работе исследуется ситуация, когда интенсивность входящего потока, а также интенсивность обслуживания зависят от  {\it полного состояния системы}, а именно -- от того времени, которое прошло с момента поступления последнего пришедшего требования и от того времени, которое имеющиеся в системе требования находились на обслуживании.

Как уже говорилось, несмотря на наличие указанной зависимости случайных величин (времена между поступлениями требований и времена обслуживания требований), мы предполагаем некоторые ограничения для интенсивности поступления требований в систему и интенсивностей их обслуживания.

Предложенный выше вариант зависимости времён обслуживания и входящего потока позволяют утверждать, что процесс обслуживания такой системы -- марковский и {\it регенерирующий}. Точками регенерации такого процесса являются моменты времени, когда все входящие в  {\it полное состояние системы} величины обнуляются.
Это может произойти тогда, когда в свободную систему приходит требование, т.е. полное состояние системы содержит нулевое время обслуживания единственного требования, и нулевое время ожидания поступления следующего требования.

\section{Описание изучаемой СМО}\label{zv16} 
Итак, мы рассматриваем систему массового обслуживания, в которой интенсивность входящего потока и интенсивность обслуживания зависят от {\it полного состояния системы}, описание которого будет дано ниже.

Это {\it полное состояние системы} в момент времени  $t$ включает: время $x^{(0)}_t$, прошедшее с момента прихода последней заявки, и времена $x^{(i)}_t$, прошедшие с момента начала обслуживания $i$-го требования из общего количества находящихся (требования нумеруются в соответствии с порядком их поступления); также для удобства мы будем использовать переменную  $n_t$ -- число заявок в системе в момент времени $t$.
Таким образом, {\it полное состояние системы} описывается вектором  $X_t=\left(n_t, x^{(0)}_t; x^{(1)}_t,x^{(2)}_t,\ldots,x^{(n_t)}_t\right)$ -- поскольку $i$-я заявка поступила ранее $(i+1)$-й, ясно, что $x^{(i)}_t\geqslant  x^{(i+1)}_t$ для $i=1,\ldots,n_t$.
Легко видеть, что  $x^{(0)}_t\leqslant  x^{(n_t)}_t$, то есть заявка с номером  $n_t$ появилась в системе не позднее времени  $t-x^{(0)}_t$.

Мы предполагаем, что интенсивность входящего потока есть функция от  $X_t$:\; $\lambda=\lambda(X_t)$ -- так же, как и интенсивность обслуживания  $i$-й заявки (которая уже обслуживалась в течении времени $x^{(i)}_t$) -- это функция $h_i=h_i(X_t)$.
Это значит, что для малого промежутка времени $\Delta>0$ при заданном $X_t=\left(n_t, x^{(0)}_t; x^{(1)}_t,x^{(2)}_t,\ldots,x^{(n_t)}_t\right)$ -- напомним, что $x^{(i)}_t\geqslant  x^{(i+1)}_t$, -- верно следующее:
\begin{equation}\label{zv8} 
\left.
\begin{array}{l}
\mathsf{P} \left\{
\begin{array}{l}
n_{t+\Delta}=n_t+1,\; x_{t+\Delta}^{(0)}=  x^{(n_{t+\Delta})}_{t+\Delta}\in(0;\Delta),
\\
\mbox{и $x_{t+\Delta}^{(i)}= x_{t}^{(i)}+\Delta$ для всех $i=1,\ldots,n_t$ }
\\
x^{(n_t+1)}_{t+\Delta}\in (0;\Delta)
\end{array}
\right\} =\lambda(X_t)\Delta+o(\Delta);
 \\
\mathsf{P} \{n_{t+\Delta}=n_t-1\}=\displaystyle  \sum\limits _{i=1}^{n_t}h_i (X_t)\Delta+o(\Delta);
\\ 
\mathsf{P} \left\{
\begin{array}{l}
 n_{t+\Delta}=n_t-1, \mbox{ и}:
 \\
 \mbox{1. для всех $j<i$, } x_{t+\Delta}^{(j)}=x_{t}^{(j)}+\Delta, 
 \\ 
 \mbox{2. для всех $j>i$ } x_{t+\Delta}^{(j)}=x_{t}^{(j+1)}+\Delta
\end{array}
\right\}=h_i(X_t)\Delta+o(\Delta)
\\
\hspace{5cm}\mbox{ -- для всех $i=1,2,\ldots,n_t$};
\\ 
\mathsf{P} \left\{n_{t+\Delta}=n_t; x_{t+\Delta}^{(i)}=x_{t}^{(i)}+\Delta \mbox{ для всех $i=0,\ldots,n_t$}\right\}=
\\
\hspace{4cm}=1-\left(\lambda(X_t)+\displaystyle  \sum\limits _{i=1}^{n_t}h_i (X_t)\right)\Delta+o(\Delta).
\end{array}
\right\}
\end{equation}
Из  (\ref{zv8}) следует, что  распределение $X_{t+\Delta}$ зависит только отраспределения $X_t$, т.е.   $X_t$ --  марковский процесс на пространстве состояний
$\displaystyle
\mathscr{X} \stackrel{\text{\rm def}}{=\!=} \bigcup\limits_{i=0}^\infty\mathscr{S} _i,
$
где $\mathscr{S} _i\stackrel{\text{\rm def}}{=\!=}  \{i\}\times\prod\limits_{j=0}^i \mathbb{R}_+$, $i\in\mathbb{Z}_+$; множество $\mathscr{S} _0$  -- это множество свободных состояний СМО (когда заявки в системе отсутствуют).

Обозначаем $\mathscr{P} _t^{X_0}$ -- распределение процесса $X_t$ с начальным состоянием $X_0$: $\mathscr{P} _t^{x}(S)=\mathsf{P} \{X_t\in S|X_0=x\}$, $x\in \mathscr{X} $, $S\in\mathscr{B}( \mathscr{X})$.

Как уже говорилось, моменты попадания процесса $X_t$ в состояние $(1,0;0)$, т.е. выхода процесса $X_t$ из множества $\mathscr{S} _0$ (в свободную систему поступило требование)  являются моментами регенерации процесса $X_t$.
Поэтому, если время между моментами регенерации имеет конечное математическое ожидание, то процесс $X_t$ -- эргодический, и его распределение слабо сходится к некоторому стационарному (инвариантному) распределению.

\begin{remark}
Предложенным способом может быть описан широкий класс СМО, например, рассмотренные выше СМО $M|M|\infty$, $M|G|\infty$, $M|G|m|0$, и другие.\hfill \ensuremath{\triangleright}
\end{remark}

Для того, чтобы среднее периода регенерации хорошо изученной СМО $M|G|\infty$ было конечным, достаточно, чтобы среднее времени обслуживания было конечным.
При этом скорость сходимости величин $\mathscr{P} _k(t)$ к их стационарным значениям $\mathscr{P} _k$ тем выше, чем больше конечных моментов у распределения времени обслуживания -- см., например, \cite{Zv3}.

Естественно поэтому рассматривать ситуацию, когда (зависимые) времена обслуживания имеют конечные моменты, а характеристики входящего потока могут быть сравнены с характеристиками пуассоновского потока.

Наложим на интенсивности следующие условия:
\begin{equation}\label{zv10} 
\begin{array}{l}
h_i(X_t)\geqslant \displaystyle \frac{K}{1+x^{(i)}_t},\;\; K>2;\quad 
0<\lambda_0\leqslant\lambda(X_t)\leqslant \Lambda<\infty.
\end{array}
\end{equation}
Из первого условия следует, что ф.р. времени обслуживания $G_i(s)$ удовлетворяет неравенству
$
\displaystyle G_i(s)\geqslant 1-\displaystyle  \frac{1}{(1+s)^K},
$
т.е. время обслуживания имеет $k\in[0;K)$ конечных моментов (см. Замечание \ref{zv9}).

Также плотность ф.р. интервалов времени поступлениями заявок $F(s)$ удовлетворяет неравенству
$
\displaystyle f_i(s)\leqslant \Lambda e^{-\lambda_0 s},
$
и распределение времён $\eta_i$ между поступлениями заявок имеет любые степенн\'ые моменты:
$
\displaystyle  \mathbb{E} \,(\eta_i)^k\leqslant\displaystyle  \frac{k!\Lambda}{\lambda_0^{k+1}},
$
в частности, эта оценка верна для {\it свободных периодов} $X_t$, т.е. промежутков времени, где $n_t\equiv 0$.\hfill \ensuremath{\triangleright}

Наша цель -- найти вычисляемые оценки скорости сходимости распределения процесса $X_t$ к стационарному распределению в условиях (\ref{zv10}).

Для этого мы воспользуемся {\it методом склеивания} -- см. Раздел \ref{zv11}.

Однако сначала мы напомним некоторые сведения из теории массового обслуживания и получим необходимые в дальнейшем оценки.

\section{Некоторые сведения о системе $MI|GI|\infty$}\label{zv20} 
Здесь мы приведём известные результаты о системе $MI|GI|\infty$ (см. \cite{Zv9}, \cite{Zv10} и др.).

В систему $MI|GI|\infty$ поступает пуассоновский поток заявок с постоянной интенсивностью $\Lambda$.
Пришедшая заявка сразу начинает обслуживание.
Её время обслуживания является случайной величиной (обозначаемой $\xi_i$) с ф.р. $G(s)$.
Все случайные величины (времена обслуживания и времена между поступлениями заявок) независимы в совокупности.

\textit{В начальный момент времени СМО свободна.}

Предполагается, что существует некоторое количество конечных моментов сл.в. $X_i$:
$
\displaystyle \mathbb{E} \,(\xi_i)^k=k\int\limits_0^\infty s^{k-1}(1-G(s))\,\mathrm{d}\, s=m_k<\infty; \;\; K\geqslant 1.
$
В дальнейшем положим
\begin{equation}\label{zv21} 
G(s)=1-\displaystyle \frac{1}{(1+s)^{K}},\qquad K>2.
\end{equation}

Условие (\ref{zv21}) гарантирует наличие $K-\varepsilon$ ($\varepsilon>0$) конечных моментов; очевидно, что
$
\displaystyle \mathbb{E} \,\xi_i=\frac{1}{K-1};\quad \mathbb{E} \,(\xi_i)^2=\frac{2}{K^2-3K+2};\mbox{ и т.д.}
$

Нагрузкой СМО назовём величину $\rho\stackrel{\text{\rm def}}{=\!=} \Lambda m_1$.

Опять описываем поведение процесса обслуживания марковским процессом  $\stackrel{\circ}{X}_t=\left(\stackrel{\circ}n_t, \stackrel{\circ}x^{(0)}_t; \stackrel{\circ}x^{(1)}_t,\stackrel{\circ}x^{(2)}_t,\ldots,\stackrel{\circ}x^{(\stackrel{\circ}n_t)}_t\right)$, где $\stackrel{\circ}n_t$ -- число требований в СМО в момент времени $t$, а  $\stackrel{\circ}x^{(0)}_t$ -- время, прошедшее с момента поступления последнего требования; $\stackrel{\circ}x^{(i)}_t$ --  время, в течении которого $i$-я по порядку поступления заявки из находящихся в системе ($i=1,2,\ldots, \stackrel{\circ}n_t$) находилась на обслуживании. (Здесь и везде в дальнейшем знак $\circ$ над случайной величиной,  процессом или их числовыми характеристиками обозначает, что речь идёт о ``классической'' СМО $MI|GI|\infty$.)

Распределение процесса $\stackrel{\circ}{X}_t$ обозначаем $\stackrel{\circ}{\mathscr{P}}_t$:
$
\stackrel{\circ}{\mathscr{P}} _t(S)=\mathsf{P} \{\stackrel{\circ}{X}_t\in S\}
$
для всех $S\in\mathscr{B} (\mathscr{X} )$, и снова $\mathscr{X} \stackrel{\text{\rm def}}{=\!=}  \bigcup\limits_{i=0}^\infty\mathscr{S} _i,
$ и  $\mathscr{S} _i\stackrel{\text{\rm def}}{=\!=}  \{i\}\times \prod\limits_{j=0}^i \mathbb{R}_+$, $i\in\mathbb{Z}_+$.

В этих предположениях известно (см., например, \cite{Zv9}), что
$$
\displaystyle\mathcal P_k(t)\stackrel{\text{\rm def}}{=\!=}  \mathsf{P} \{\stackrel{\circ}n_t=k\}=\stackrel{\circ}{\mathscr{P}} _t(\mathscr{S} _0)=\displaystyle  \frac{e^{-\mathfrak{G} (t)}\big(\mathfrak{G} (t)\big)^k}{k!},
$$ где $\displaystyle\mathfrak{G} (t)\stackrel{\text{\rm def}}{=\!=} \displaystyle  \Lambda\int\limits_0^t(1-G(s))\,\mathrm{d}\, s.
$
\begin{remark}\label{zv26} 
$\mathfrak{G} (+\infty)=\Lambda m_1$, и $\mathcal P_0(t)\geqslant e^{-\Lambda m_1}=e^{-\rho}$ при $t\geqslant 0$.\hfill \ensuremath{\triangleright}
\end{remark}

Кроме этого, известно распределение периода занятости СМО  $M|G|\infty $ (см., например, \cite{Zv10}).
Обозначим $\stackrel{\circ}\zeta_i$ -- $i$-й период занятости СМО  $M|G|\infty $.
Эта сл.в. имеет ф.р.
$
\displaystyle B(x)\stackrel{\text{\rm def}}{=\!=}  \mathsf{P} \{\stackrel{\circ}\zeta_i\leqslant x\}=1- \frac{1}{\Lambda}\displaystyle  \sum\limits _{k=1}^\infty c^{n\ast}(x),
$
где  $c(x)\stackrel{\text{\rm def}}{=\!=}  \displaystyle  \Lambda(1-G(x))e^{-\mathfrak{G} (x)}$, а $c^{n\ast}$ -- это $n$-я свёртка функций $c(x)$.

Также можно выписать преобразование Лапласа для $B(x)$ (см., например, \cite{Zv10}):
\begin{equation}\label{zv17} 
  \mathcal{L}[B](s)=1+\frac{s}{\Lambda}-\frac{1}{\Lambda \displaystyle \int\limits_0^\infty \exp\left(-st-\Lambda\int\limits_0^t [1-G(v)]\,\mathrm{d}\, v\right)\,\mathrm{d}\, t}.
\end{equation}
Впрочем, моменты периода занятости можно вычислять и по формулам
\begin{equation}\label{zv18} 
\mathbb{E} \,(\stackrel{\circ}\zeta_i)^n=(-1)^{n+1}\left\{\frac{e^\rho}{\Lambda}nC^{(n-1)} - e^\rho \displaystyle  \sum\limits _{k=1}^{n-1}C_n^k \mathbb{E} \,(\stackrel{\circ}\zeta_i)^{n-k} C^{(k)} \right\},
\qquad
n\in\mathbb{N},
\end{equation}
где
$\displaystyle
C^{(k)}=\int\limits_0^\infty (-t)^k \left(\exp\left(-\Lambda
\int\limits_0^\infty [1-G_0(v)]\,\mathrm{d}\, v\right)\right)\Lambda
[1-G_0(t)]\,\mathrm{d}\, t.
$

Используя (\ref{zv17}) и (\ref{zv18}), можно получить значения
\begin{equation}\label{zv36} 
\begin{array}{l}
\hspace{2cm}\displaystyle \mathbb{E}  \, \stackrel{\circ}\zeta_i= \frac{e^\rho-1}{\Lambda} 
\mbox{ и }
 \mathbb{E}  \,(\stackrel{\circ}\zeta_i)^2=
 \\ =\displaystyle \frac{2e^ {2\rho}}{\Lambda}\displaystyle \int\limits_0^\infty \left(\exp\left(-\Lambda \int\limits_0^\infty [1-G(v)]\,\mathrm{d}\, v\right)\right)\,\mathrm{d}\, t+ \displaystyle \frac{2e^\rho} {\Lambda^2}-\frac{2}{\Lambda^2}-2\displaystyle \frac{e^\rho-1}{\Lambda^2}.
\end{array}
\end{equation}

Но для более высоких моментов  $\zeta_i$ вычисления по формулам (\ref{zv17}) и (\ref{zv18}) очень сложны и трудно применимы на практике.

Однако нам в дальнейшем потребуются оценки моментов  $\zeta_i$.
Для этого используем
\begin{proposition}\label{zv19} 
Если $\mathbb{E} \,\xi_i^k<\infty$ при некотором $k\in\mathbb{N} $, то верна оценка
\begin{equation}\label{zv25} 
\mathbb{E} \,(\stackrel{\circ}\zeta_i)^k\leqslant \frac{\mathbb{E} \,\xi_i ^k}{1-e^{-\rho}} \times \varphi((1-e^{-\rho}),k-1),
\end{equation}
где $\displaystyle  \varphi(x,k)\stackrel{\text{\rm def}}{=\!=}  \displaystyle  \sum\limits _{n=1}^\infty n^{k-1}x^n=\left(x\,\frac{\,\mathrm{d}\,}{\,\mathrm{d}\, x}\right)^k\frac{1}{1-x}$.
\hfill \ensuremath{\triangleright}
\end{proposition}

\noindent\textbf{Доказательство.}
Функция  $c(x)$, используемая в (\ref{zv17}), неотрицательна, и
\begin{multline*}
\displaystyle  \int\limits_0 ^\infty c(x)\,\mathrm{d}\, x=\int\limits_0 ^\infty \Lambda(1-G(x))\exp\left(\displaystyle  -\Lambda\int\limits_0^x(1-G(s))\,\mathrm{d}\, s\right) \,\mathrm{d}\, x =
\\
=\left(1-\exp\left(\displaystyle  -\Lambda\int\limits_0^\infty(1-G(s))\,\mathrm{d}\, s\right)\right)= \left(1-e^{\displaystyle -\Lambda m_1}\right)=(1-e^{-\rho})\stackrel{\text{\rm def}}{=\!=}  \varrho,
\end{multline*}
где опять $\rho=\displaystyle  \Lambda m_1$.

Т.е.  функция $\varsigma(x)\stackrel{\text{\rm def}}{=\!=}  \varrho^{-1}\, c(x)$ является плотностью распределения некой сл.в.  $\vartheta$, а $\varsigma^{n\ast}(x)$ -- это плотность распределения $\displaystyle \displaystyle  \sum\limits _{j=1}^n\vartheta_j$ -- суммы $n$ одинаково распределённых и независимых сл.в. с тем же распределением, что и $\vartheta$.

Поэтому $\displaystyle \int\limits_0^  \infty \varsigma^{n\ast}(x)\,\mathrm{d}\, x=1$.
Следовательно,
\begin{multline*}
  \mathbb{E} \,\stackrel{\circ}\zeta_i= \int\limits_0^\infty (1-B(x))\,\mathrm{d}\, x=\displaystyle  \sum\limits _{n=1}^\infty \frac{1}{\Lambda}\int\limits_0^\infty c^{n\ast}(x)\,\mathrm{d}\, x =
  \\
  = \frac{1}{\Lambda}\displaystyle  \sum\limits _{n=1}^\infty \varrho^n\int\limits_0^\infty \varsigma^{n\ast}(x)\,\mathrm{d}\, x=\frac{1}{\Lambda}\displaystyle  \sum\limits _{n=1}^\infty \varrho^n   =\frac{1}{\Lambda}\times \frac{\varrho}{1-\varrho}=\frac{e^{\rho}-1}{\Lambda}.
\end{multline*}

Далее, если  $m_{k+1}=\mathbb{E} \,\xi_i^{k+1}<\infty$ при некотором $k\in\mathbb{N} $, то
\begin{multline*}
\mathbb{E} \,\vartheta^k=\int\limits_0 ^\infty x^k\varsigma(x)\,\mathrm{d}\, x =\int\limits_0 ^\infty x^k \varrho^{-1} \Lambda(1-G(x))e^{-\mathfrak{G} (x)}\,\mathrm{d}\, x \leqslant
\\
\leqslant \frac{\Lambda}{\varrho }\int\limits_0^ \infty x^k (1-G(x))\,\mathrm{d}\, x =\frac{\Lambda\mathbb{E} \,\xi_i^{k+1}}{\varrho(k+1)},
\end{multline*}
и
$ \displaystyle
\int\limits_0 ^\infty x^k\varsigma^{n\ast}(x)\,\mathrm{d}\, x=\mathbb{E} \left(\displaystyle  \sum\limits _{j=1}^n\vartheta_j\right)^k\leqslant n^{k}\mathbb{E} \,\vartheta^k\leqslant \frac{ n^{k}\Lambda m_{k+1}}{\varrho (k+1)}.
$

В последней оценке мы использовали неравенство Йенсена (для $k\geqslant 1$) в таком виде:
\begin{equation}\label{zv33} 
  \displaystyle  (a_1+\ldots+a_n)^k\leqslant n^{k-1}(a_1^k+\ldots +a_n^k).
\end{equation}

Теперь оценим $\mathbb{E} \,(\stackrel{\circ}\zeta_i)^k$:
\begin{multline*}
\mathbb{E} \,(\stackrel{\circ}\zeta_i)^k=\int\limits_0 ^ \infty k x^{k-1}(1-B(x))\,\mathrm{d}\, x
=\frac{1}{\Lambda}\displaystyle  \sum\limits _{n=1}^\infty\int\limits_0^ \infty  k x^{k-1} c^{n\ast}(x)\,\mathrm{d}\, x =
\\
= \frac{1}{\Lambda}\displaystyle  \sum\limits _{n=1}^\infty\int\limits_0^ \infty  k x^{k-1} \varrho^n\varsigma^{n\ast}(x)\,\mathrm{d}\, x= \frac{1}{\Lambda}\displaystyle  \sum\limits _{n=1}^\infty  \varrho^n k\, \mathbb{E} \left(\displaystyle  \sum\limits _{j=1}^n \vartheta_j\right)^ {k-1}\leqslant
\\
\leqslant \frac{1}{ \Lambda} \displaystyle  \sum\limits _{n=1}^\infty \varrho^nk n^{k-1}\frac{\Lambda\mathbb{E} \, \xi_i^k}{\varrho k}= \frac{\mathbb{E} \,\xi_i ^k}{\varrho } \displaystyle  \sum\limits _{n=1}^\infty n^{k-1} \varrho^n= \frac{\mathbb{E} \,\xi_i ^k}{\varrho } \times \varphi(\varrho,k-1).
\end{multline*}
Предложение \ref{zv19} доказано.
\hfill \ensuremath{\blacksquare}

\section{Оценки для процесса $X_t$}\label{zv31} 
Вернёмся к исследуемому нами процессу $X_t$, описанному в Разделе \ref{zv16}.
Предполагаем, что в начальный момент времени $X_0=(0,0)$.
Т.е. система свободна и начинается отсчёт времени до поступления первого требования. В момент поступления первоготребования процесс $X_t$ попадает в состояние $(1,0;0)$ -- это момент регенерации.

Напомним, что описанный в разделе \ref{zv20} процесс $\stackrel{\circ}{X}_t$ также стартует из этого состояния: $\stackrel{\circ}X_0=(0,0)$. В момент поступления первого требования начинается первый период регенерации.
\begin{definition}
Случайная величина $\eta$ не превосходит случайную величину $\theta$ \textit{по распределению}, если при всех $s\in\mathbb{R}$ верно неравенство
$ \displaystyle
F_\eta(s)=\mathsf{P} \{\eta\leqslant s\} \geqslant \mathsf{P} \{\theta\leqslant s\}=F_\theta(s),
$
откуда следует, что ф.р. сл.в. $\theta$ не превосходит ф.р. сл.в. $\eta$.
Это -- отношение порядка, будем обозначать его $\eta\prec \theta$.\hfill \ensuremath{\triangleright}
\end{definition}
\begin{proposition}\label{zv22} 
$\eta\prec \theta$ тогда и только тогда, когда на некотором вероятностном пространстве $(\Omega,\mathscr{F} ,\mathsf{P} )$ можно построить сл.в. $\eta'$ и $\theta'$ такие, что $\eta'\leqslant\theta'$ при всех $\omega\in\Omega$, и $\mathsf{P} \{\eta\leqslant s\}=\mathsf{P} \{\eta'\leqslant s\}$ и $\mathsf{P} \{\theta\leqslant s\}=\mathsf{P} \{\theta'\leqslant s\}$ для всех $s\in\mathbb{R}$.
\hfill \ensuremath{\triangleright}
\end{proposition}

\noindent\textbf{Доказательство.}
  См. \cite{sht}.
\hfill \ensuremath{\blacksquare}
\begin{proposition}\label{zv23} 
Если $X_0=\stackrel{\circ}X_0=(0,0)$, и выполнено условие (\ref{zv10})  (напомним, что входящий поток для $\stackrel{\circ}{X}_t$ -- пуассоновский и имеет интенсивность $\Lambda$, а время обслуживания имеет ф.р. (\ref{zv21})), то при всех $t\geqslant 0$ выполнено соотношение:
\begin{equation}\label{zv24} 
n_t\prec \stackrel{\circ}n_t.
\end{equation}
\hfill \ensuremath{\triangleright}
\end{proposition}

\noindent\textbf{Доказательство.}
Для доказательства используется \textit{метод одного вероятностного пространства}, который заключается в том, что на одном вероятностном пространстве строятся процессы $X_t$ и $\stackrel{\circ}X_t$ (а точнее говоря, их версии) таким образом, что на всех траекториях процессов выполнено неравенство $n_t\leqslant \stackrel{\circ}n_t$ -- при всех $t\geqslant 0$.

Эта конструкция возможна по двум причинам.

{\bf I.} Входящий поток для процесса $X_t$ имеет интенсивность, не превосходящую интенсивность входящего потока для $\stackrel{\circ}X_t$ ($\lambda(X_t)\leqslant \Lambda$ из условия (\ref{zv10})).

Поэтому можно строить на одном вероятностном пространстве два входящих потока для обоих процессов: поток $\mathfrak{F}_t$ с интенсивностью $\lambda(X_t)$ для процесса $X_t$, и поток $\mathfrak{F}^+_t$ с интенсивностью $\Lambda-\lambda(X_t)$; объединение потоков $\mathfrak{F}_t$  и $\mathfrak{F}^+_t$ -- это поток с интенсивностью $\Lambda$ для процесса $\stackrel{\circ}X_t$.

Таким образом, в каждый момент времени $t$ число заявок, поступивших в СМО с процессом $X_t$ не превосходит число заявок, поступивших в СМО с процессом $\stackrel{\circ}X_t$.
При этом в моменты поступления заявок в в СМО с процессом $X_t$ одновременно поступают заявки в СМО с процессом $\stackrel{\circ}X_t$; назовём такие заявки ``общими''.

Иначе говоря, в каждый момент времени вектор $\left( x^{(1)}_t,x^{(2)}_t,\ldots,x^{(n_t)}_t\right)$ состоит из компонент вектора $\left( \stackrel{\circ}x^{(1)}_t,\stackrel{\circ}x^{(2)}_t,\ldots,\stackrel{\circ}x^{(\stackrel{\circ}n_t)}_t\right)$, т.е. $\left( x^{(1)}_t,x^{(2)}_t,\ldots,x^{(n_t)}_t\right)$ состоит из обслуженных времён находящихся на обслуживании требований, совпадающих с некоторыми из элементов вектора $\left( \stackrel{\circ}x^{(1)}_t,\stackrel {\circ}x^{(2)}_t,\ldots,\stackrel {\circ}x^{(\stackrel{\circ}n_t)}_t \right)$.

{\bf II.} Интенсивность обслуживания $i$-го требования СМО, описываемой процессом $\stackrel{\circ}X_t=\left(\stackrel{\circ}n_t, \stackrel{\circ}x^{(0)}_t; \stackrel{ \circ}x^{(1)}_t, \stackrel{\circ}x^{(2)}_t,\ldots, \stackrel{\circ}x^{(\stackrel{ \circ}n_t)}_t\right)$, равна
$ \displaystyle
\stackrel{\circ}h_i(\stackrel {\circ}X_t)= \frac{G'\left(\stackrel{ \circ}x^{(i)}_t\right)} {1-G\left(\stackrel{ \circ}x^{(i)}_t\right)} =\frac{K}{1+\stackrel{ \circ}x^{(i)}_t}.
$
Таким образом,  интенсивность обслуживания ``общих'' заявок у процесса $\stackrel{\circ}X_t$ не превосходит интенсивность обслуживания этих же заявок у процесса $X_t$, и это -- второе обстоятельство, позволяющее применить метод одного вероятностного пространства.

Учитывая, что время обслуживания имеет ф.р., которую можно восстановить по интенсивности обслуживания (см. (\ref{zv2})--(\ref{zv3})), получаем, что остаточное время обслуживания ``общих'' заявок у процесса $X_t$ не превосходит\textit{ по распределению }остаточное время обслуживания этих заявок у процесса $\stackrel{\circ}X_t$.
Поэтому можно построить на одном вероятностном пространстве эти случайные величины так, что ``общие'' заявки завершат обслуживание у процесса $X_t$ не позднее, чем у процесса $\stackrel{\circ}X_t$ -- см. Предложение \ref{zv22}.

Стало быть, в СМО с процессом $\stackrel{\circ}X_t$ поступает больший поток заявок, чем в СМО с процессом $X_t$, и одновременно поступившие (общие) заявки покидают СМО с процессом $X_t$ раньше.
Из этого следует утверждение предложения \ref{zv23}.
\hfill \ensuremath{\blacksquare}
\begin{corollary}
Из (\ref{zv24}) следует, что свойство упорядоченности по распределению верно не только для $n_t$ и $\stackrel{\circ}n_t$, но и для периодов занятости процессов $X_t$ и  $\stackrel{\circ}X_t$.
Учитывая, что процессы $X_t$ и $\stackrel{\circ}X_t$ регенерирующие (напомним, моментом регенерации является момент, когда процесс принимает значение $(1,0;0)$, т.е. момент начала периода занятости), распределение периодов занятости $\zeta_i$ процесса $X_t$ и периодов занятости $\stackrel{\circ}\zeta_i$ процесса $\stackrel{\circ}X_t$ не зависят от номера $i$.
То есть $
\zeta_i\prec\stackrel{\circ}\zeta_t,
$
и для всех $k>0$ верно неравенство
$
\displaystyle\boxed{\mathbb{E} (\zeta_i)^k\leqslant\mathbb{E} (\stackrel{\circ}\zeta_i)^k}.
$
То есть оценки (\ref{zv25}) верны и для $\zeta_i$.\hfill \ensuremath{\triangleright}
\end{corollary}
\begin{corollary}\label{zero}  
Из Замечания \ref{zv26} и (\ref{zv24}) следует, что при всех $t\geqslant 0$ (и $X_0=(0,0)$) верно неравенство
\begin{equation}\label{zv32}
\mathsf{P} \{X_t\in\mathscr{S} _0\}=\mathsf{P} \{n_t=0\}\geqslant \mathcal{P}_0(t)\geqslant e^{-\rho}.
\end{equation}
\hfill \ensuremath{\triangleright}
\end{corollary}

\section{Метод склеивания}\label{zv11} 

Первоначально метод склеивания марковских процессов был применён к цепям Маркова с конечным числом состояний (см. \cite{Zv6,Zv14}).
Основная идея этого метода состоит в следующем.

Пусть два однородных марковских процесса $Y_t$ и $Z_t$ с одинаковой переходной функцией в начальный момент времени имеют различные состояния; они заданы на $(\Omega,\mathscr{F} ,\mathsf{P} )$, и измеримы относительно фильтрации $\mathscr{F} _t$.

Обозначим $\tau=\tau(Y_0,Z_0)\stackrel{\text{\rm def}}{=\!=}  \inf\{t>0:\, Y_t=Z_t\}$.
То есть $Y_\tau=Z_\tau$, и при $t>\tau$ распределения процессов $Y_t$ и $Z_t$ совпадают в соответствии с марковским свойством, т.е.
$
\mathsf{P} \{Y_t\in S|t>\tau\}=\mathsf{P} \{Z_t\in S|t>\tau\}
$
для всех $S\in\mathscr{B} (\mathscr{X} )$, $Y_t$ and $Z_t$, здесь $\mathscr{B} (\mathscr{X} )$ -- $\sigma$-алгебра на пространстве состояний $\mathscr{X}$ процессов $Y_t$ и $Z_t$.

Имеет место {\it основное неравенстно склеивания}:
\begin{multline}\label{zv12}  
|\mathsf{P} \{Y_t\in S\} - \mathsf{P} \{Z_t\in S\}|=
\\
=|\mathsf{P} \{Y_t\in S\;\&\; \tau>t\} - \mathsf{P} \{Z_t\in S\;\&\; \tau>t\} +
\\
+ |\mathsf{P} \{Y_t\in S\;\&\; \tau<t\} - \mathsf{P} \{Z_t\in S\;\&\; \tau<t\}|=
\\
=|\mathsf{P} \{Y_t\in S\;\&\; \tau>t\} - \mathsf{P} \{Z_t\in S\;\&\; \tau>t\} |\leqslant \mathsf{P} \{\tau>t\},
\end{multline}
поскольку при $\tau<t $ ``склеивание'' процессов $Y_t$ и $Z_t$ уже произошло, их распределения в момент $t$ одинаковы, и $\mathsf{P} \{Y_t\in S| \tau<t\}=\mathsf{P} \{Z_t\in S| \tau<t\}$; здесь $S\in\mathscr{B} (\mathscr{X} )$.

Поэтому если для некоторой возрастающей положительной функции $\phi(t)$ удастся найти оценку $\mathbb{E} \, \phi(\tau(Y_0,Z_0))\leqslant C(X_0,Y_0)$, то можно применить неравенство Маркова:
\begin{equation}\label{zv13} 
\mathsf{P} \{\tau(Y_0,Z_0)>t\}=\mathsf{P} \{\phi(\tau(Y_0,Z_0))>\phi(t)\}\leqslant \frac{\mathbb{E} \, \phi(\tau(Y_0,Z_0))}{\phi(t)}\leqslant \frac{C(X_0,Y_0)}{\phi(t)}.
\end{equation}

Далее, если существует стационарное инвариантное распределение $\mathscr{P} $, к которому (слабо) сходятся распределения  $Y_t$ и $Z_t$ при $t\to\infty$, то из (\ref{zv12}) и (\ref{zv13}) получаем
\begin{equation}\label{zv35} 
\left.\begin{array}{l}
\hspace{1cm}  |\mathsf{P} \{Y_t\in S\} - \mathscr{P} (S)|\leqslant \displaystyle \int\limits_\mathcal{\mathscr{X} } \frac{C(X_0,y)}{\phi(t)}\mathscr{P} (\,\mathrm{d}\, y)=\displaystyle \frac{\hat C(X_0)}{\phi(t)},
\\
\mbox{откуда}
\\
\hspace{1cm} \|\mathscr{P} _t^{X_0}-\mathscr{P} \|_{TV}\stackrel{\text{\rm def}}{=\!=}   \sup\limits _{S\in \mathscr{B} (\mathscr{X} )}|\mathscr{P} _t^{X_0}(S)-\mathscr{P} (S)|\leqslant 2 \displaystyle \frac{\hat  C(X_0)}{\phi(t)}.
\end{array}
\right\}
\end{equation}

Изначально метод склеивания был предложен для цепей Маркова с дискретным множеством состояний (\cite{Zv6}), однако изучаемый нами процесс $X_t$ с переходными вероятностями (\ref{zv8}) -- это процесс в непрерывном времени, и для двух независимых процессов $X_t$ и $X_t'$, удовлетворяющих (\ref{zv8}), но с различными начальными состояниями $\mathsf{P} \{\tau(X_0,x_0')<\infty\}=0$.

Поэтому мы будем использовать понятие \emph{успешной склейки}, предложенное в \cite{Zv5}.

\begin{definition} \label{zv30} 
Пусть $X_t$ и $X_t'$ -- две независимые версии одного марковского процесса с пространством состояний $\mathscr{X} $: $X_t$ с начальным состоянием $X_0$, и $X_t'$ с начальным состоянием $X_0'$.

\textsc {Успешной склейкой } (см. \cite{Zv5}) процессов $X_t$ и $X_t'$ называется парный процесс  $\mathcal{Z}_t=(Y_t,Y_t')$, заданный на некотором вероятностном пространстве, такой, что
\begin{itemize}
  \item [(i)] Для всех $s\geqslant 0$ и  $S\in \mathscr{B} (\mathscr{X} )$ выполнены равенства

  $
   \mathbf{P} \{Y_t\in S\}= \mathbf{P} \{X_t\in S\}$, $ \mathbf{P} \{Y'_t\in S\}= \mathbf{P} \{X'_t\in S\}.
  $
  (Соответственно, $Y_0=X_0$ и $Y'_0=X'_0$.)

{При этом для любого фиксированного момента $t\geqslant 0$ случайные величины $Y_t$, $Y_t'$, $X_t$ и $X_t'$ рассматриваются как некоторые случайные величины, но не как случайные процессы! Т.е. речь идёт не о совпадении распределений двух процессов, под которым понимается совпадение всех конечномерных распределений.
}
\item[(ii)] Для всех $t>\tau=\tau(X_0,X'_0)=\tau(Y_0,Y'_0) \stackrel{\text{\rm def}}{=\!\!=\!\!=}\inf\{t\geqslant 0:\,Y_t=Y_t'\}$ выполняется равенство $Y_t=Y_t'$.
\item[(iii)] Для любых начальных значений $X_0,X_0'\in  \mathscr{X}$ верно: $ \mathbf{P} \{\tau(X_0,X'_0)<\infty\}=1$.
\hfill \ensuremath{\triangleright}
\end{itemize}
\end{definition}
Для  \emph{успешной склейки} можно также использовать \textit{основное неравенство склеивания} (\ref{zv12}):
\begin{multline*}
  |\mathsf{P} \{X_t\in S\}-\mathsf{P} \{X_t'\in S\}|=|\mathsf{P} \{Y_t\in S\}-\mathsf{P} \{Y_t'\in S\}|
\\
=|\mathsf{P} \{Y_t\in S\}-\mathsf{P} \{Y_t'\in S\}|\times (\mathbf{1} (\tau> t)+ \mathbf{1} (\tau\leqslant t))\leqslant \mathsf{P} \{\tau> t\}.
\end{multline*}

Итак, наша цель -- сконструировать успешную склейку для двух версий процесса $X_t$ ($X_t$ с начальным состоянием $X_0=(0,0)$ и $X_t'$ с начальным состоянием $X_0'=\left(n_0', x_t'{}^{(0)}; x_t'{}^{(1)},x_t'{}^{(2)},\ldots,x_t'{}^{(n'_0)}\right)$), причём каждый из процессов удовлетворяет условиям (\ref{zv8}), и из этого получить оценку скорости сходимости распределения  процесса $X_t$ к стационарному распределению.

Для применения метода склеивания нам понадобится следующая
\subsection{Основная лемма склеивания}
Этот факт изложен во многих публикациях (см., например, \cite{Zv13,verbut}), но здесь мы приведём его простейшую форму.

\begin{definition} Общей частью распределения двух случайных величин  $\xi_1$ и $\xi_2$ с функциями распределения  $\Psi_j(s)=\mathbf{P}\{\xi_j\leqslant s\}$  $s\in \mathbb{R}$ ($j=1,2$), назовём\linebreak
$\displaystyle
\varkappa \stackrel{{\rm def}}{=\!\!\!=} \varkappa(\Psi_1(s),\Psi_2(s))  \stackrel{{\rm def}}{=\!\!\!=}  \displaystyle \int\limits _{-\infty}^\infty \min(\psi_1(u),\psi_2(u))\,\mathrm{d}\, u,
$
где\\

$\psi_j(s)  \stackrel{{\rm def}}{=\!\!\!=}
\left\{
\begin{array}{ll}
 \Psi_j'(s), &\mbox{ если существует }\Psi_j'(s),
 \\
 0, &  \mbox{ в противном случае, }
\end{array}
\right.
$
\;\; здесь   $j=1,2$.\hfill \ensuremath{\triangleright}
\end{definition}
\begin{lemma}\label{zv29} 
Если  $\varkappa>0$, то {\it на некотором вероятностном пространстве} можно построить две случайные величины $\xi_1'$ и $\xi_2'$ такие, что:
\begin{equation}\label{zv14}
\mathsf{P} \{\xi_j'\leqslant s\}=\mathsf{P} \{\xi_j\leqslant s\}=\Psi_j;
\end{equation}
\begin{equation}\label{zv15}
\mathsf{P} \{\xi_1'=\xi_2'\}\geqslant \varkappa.
\end{equation}
\hfill \ensuremath{\triangleright}
\end{lemma}

\textbf{Доказательство.}
Обозначим
$$\begin{array}{ll}
\psi(u)\stackrel{\text{\rm def}}{=\!=}  \min(\psi_1(u),\psi_2(u));&
\Psi(s)\stackrel{\text{\rm def}}{=\!=} \displaystyle  \int\limits_{-\infty}^s \psi(u)\,\mathrm{d}\, u;
\\
\hat{\Psi}(s)\stackrel{\text{\rm def}}{=\!=}  \displaystyle \frac{\Psi(s)}{\varkappa};
&
\hat{\Psi}_j(s)\stackrel{\text{\rm def}}{=\!=}  \displaystyle \frac{\Psi_j(s)-\Psi(s)}{1-\varkappa}, \; j=1,2
\end{array}
$$
(если $\varkappa=1$, то $\hat{\Psi}_j(s)\equiv 0$).

Пусть $\mathscr{U} _1$, $\mathscr{U} _2$ и  $\mathscr{U} _3$ -- независимые случайные величины, равномерно распределённые на отрезке $[0;1]$.
Положим
$ \displaystyle
\xi_j'\stackrel{\text{\rm def}}{=\!=}  \hat{\Psi}^{-1}(\mathscr{U} _j)\times \mathbf{1} (\mathscr{U} _3\leqslant\varkappa) + \hat{\Psi}^{-1}_j(\mathscr{U} _j)\times \mathbf{1} (\mathscr{U} _3>\varkappa)
$.

Напомним, что под обратной функцией монотонной функции $f(s)$ мы понимаем
$
f^{-1}(u)\stackrel{\text{\rm def}}{=\!=} \inf\{s:\;f(s)\geqslant u\}.
$

Случайные величины $\xi_j'$ удовлетворяют условиям (\ref{zv14})--(\ref{zv15}).

Лемма доказана.
\hfill \ensuremath{\blacksquare}

\begin{remark}\label{zv27} 
Если $n_\vartheta=0$ в некоторый момент времени $\vartheta$ (СМО свободна), то остаточная длина свободного периода процесса $X_t$ (описанного в Разделе \ref{zv16}) имеет плотность распределения
$ \displaystyle
\varphi(s)=\lambda(X_{\vartheta+s})\exp\left(\displaystyle \int\limits_0^s \lambda(X_{\vartheta+u})\,\mathrm{d}\, u\right)\in(\lambda_0e^{-\Lambda s}; \Lambda_0e^{-\lambda s}).
$
Поэтому, если рассматриваются два процесса $X_t$ и $X_t'$ с различными начальными состояниями $X_0$ и $X_0'$, и в некоторый момент времени $\vartheta$ они оба оказались в множестве $\mathscr{S} _0$ (СМО свободны), то общая часть $\varkappa_\vartheta$ распределений их остаточных времён пребывания в множестве $\mathscr{S} _0$ удовлетворяет неравенству
\begin{equation}\label{zv28} 
  \widehat\varkappa_\vartheta\geqslant \int\limits_0^\infty \lambda_0e^{-\Lambda s}\,\mathrm{d}\, s=\frac{\lambda_0}{\Lambda}\stackrel{\text{\rm def}}{=\!=} \widehat\varkappa.
\end{equation}
\hfill \ensuremath{\triangleright}
\end{remark}

\section{Конструирование успешной склейки для пары процессов $X_t$, $X_t'$}
Рассмотрим два независимых процесса:  процесс $X_t$ с начальным состоянием $X_0=(0,0)$ и процесс $X_t'$ с начальным состоянием
$ \displaystyle
X_0'=\linebreak\left(n_0',x'_0{}^{(0)}; x'{}^{(1)}_0,x_0'{}^{(2)},\ldots,x_0'{}^{(n_0')}\right).
$

Будем строить их склейку $\mathscr{Z} _t=(Y_t,Y_t')$ на вероятностном пространстве $(\Omega,\mathscr{F} , \mathsf{P} )=\prod\limits_{i=0}^\infty (\Omega_i, \mathscr{F} _i , \mathsf{P} _i)$, где $\Omega_i=[0;1]$, $\mathscr{F} _i = \mathscr{B} ([0;1])$, и $\mathsf{P} _i $ -- мера Лебега на $[0;1]$.

Опишем процедуру конструирования склейки по шагам, которые производятся в моменты изменения первой компоненты любого из процессов $Y_t=\left(m_t,y_t^{(0)}; y^{(1)}_t,y_t^{(2)},\ldots,y_t^{(m_t)}\right)$, $Y_t'=\left(m_t',y_t'{}^{(0)}; y_t'{}^{(1)},y_t'{}^{(2)},\ldots,y_t'{}^{(m_t')}\right)$ -- обозначим эти моменты $t_k$:
$ \displaystyle
t_0\stackrel{\text{\rm def}}{=\!=}  0; \;\; t_k\stackrel{\text{\rm def}}{=\!=}  \inf\{t>t_{k-1}: m_{t}+m'_{t}\neq m_{t_{k-1}}+m'_{t_{k-1}} \}.
$

На каждом шаге строятся остаточные \emph{витруальные} времена обслуживания находящихся в СМО заявок (в предположении, что до окончания обслуживания данной заявки первые компоненты -- $n_t$ и $n_t'$ -- обоих процессов не изменятся) и остаточные \emph{витруальные} времена до прихода следующей заявки.
Для этого используется некоторое количество $\nu_k$ экземпляров $(\Omega_i, \mathscr{F} _i , \mathsf{P} _i)$ -- так, что на каждом шаге используются новые $(\Omega_i, \mathscr{F} _i , \mathsf{P} _i)$.

Опишем алгоритм конструирования.

\textbf{А. }Если $m_{t_k}+m'_{t_k}>0$, для удобства положим $Y_{\vartheta_k}=(\imath,\alpha_0;\alpha_1, \alpha_2, \ldots \alpha_\imath)$, $Y_{\vartheta_k}'=(\jmath,\beta_0;\beta_1, \beta_2, \ldots \beta_\jmath)$.

На $k$-том шаге строятся сл.в. на пространствах $(\Omega_i, \mathscr{F} _i , \mathsf{P} _i)$, $i=\ell_k+1, \ell_k+2,\ldots \ell_{k+1}$; при этом $\nu_k\stackrel{\text{\rm def}}{=\!=}  \ell_{k+1}-\ell_k=n_{\vartheta_k}+ n_{\vartheta_k}'+2= \imath+\jmath+2$.
К этому моменту времени было использовано $\ell_k=\displaystyle  \sum\limits _{i=0}^{k-1}\nu_i$ экземпляров пространства $(\Omega_i, \mathscr{F} _i , \mathsf{P} _i)$.

В момент $t_k$ \emph{виртуальное }время до поступления очередного требования $\xi$ (при условии, что не закончится обслуживания какого-то требования) для $Y_t$ имеет ф.р. $F_1(s)=1-\exp\left(\displaystyle \int\limits_0^s\lambda(\imath,\alpha_0+u;\alpha_1+u, \alpha_2+u, \ldots \alpha_\imath+u)\,\mathrm{d}\, u\right)$.
Строим на пространстве $(\Omega_{\ell_k+1}, \mathscr{F} _{\ell_k+1} , \mathsf{P} _{\ell_k+1})$ сл.в. $\xi\stackrel{\text{\rm def}}{=\!=}  F_1^{-1}(\omega)$, $\omega\in \Omega_{\ell_k+1}$.

В момент $t_k$ \emph{виртуальное } время до окончания обслуживания $j$-го ($j=1,2,\ldots,\imath$) требования $\eta_j$ (при условии, что не закончится обслуживания какого-то другого требования и не поступит новых требований)  $Y_t$ имеет ф.р. \\ $G_{1,j}(s)=1-\exp\left(\displaystyle \int\limits_0^s h_j(\ell,\alpha_0+u;\alpha_1+u, \alpha_2+u, \ldots \alpha_\imath+u)\,\mathrm{d}\, u\right)$.

Строим на пространстве $(\Omega_{\ell_k+j+1}, \mathscr{F} _{\ell_k+j+1} , \mathsf{P} _{\ell_k+j+1})$ сл.в. $\eta_j\stackrel{\text{\rm def}}{=\!=}  G_{1,j}^{-1}(\omega)$, $\omega\in \Omega_{\ell_k+j+1}$.

Аналогично строятся сл.в. $\xi'$ и $\eta_i'$ ($i=1,2,\ldots,\jmath)$:
\\
$F_2(s)=1-\exp\left(\displaystyle \int\limits_0^s\lambda(\jmath,\beta_0+u;\beta_1+u, \beta_2+u, \ldots \beta_\jmath+u)\,\mathrm{d}\, u\right)$, $\xi'\stackrel{\text{\rm def}}{=\!=}  F_2^{-1}(\omega)$, $\omega\in \Omega_{\ell_k+\imath +1}$;
\\
$G_{2,j}(s)=1-\exp\left(\displaystyle \int\limits_0^s h_j(\jmath,\beta_0+u;\beta_1+u, \beta_2+u, \ldots \beta_\jmath+u)\,\mathrm{d}\, u\right)$, $\eta_j'\stackrel{\text{\rm def}}{=\!=}  G_{2,j}^{-1}(\omega)$, $\omega\in \Omega_{\ell_k+\imath+j+2}$.

Использовано $\nu_k=\imath+\jmath+2$ экземпляров $(\Omega_i, \mathscr{F} _i , \mathsf{P} _i)$; $\ell_{k+1}=\ell_k+\nu_k$.

В момент времени $t_{k+1}=t_k+\min\{\xi,\xi',\eta_j,\eta'_i \; (j=1,2,\ldots \imath,\; i=1,2,\ldots \jmath)\}$ произойдёт изменение первой компоненты одного из процессов $Y_t$ и $Y_t'$.

Если опять $m_{t_k}+m'_{t_k}>0$, повторяем процедуру \textbf{А}.

\textbf{Б.} Если $m_{t_k}+m'_{t_k}=0$, то $Y_{t_k}=(0,\alpha_0)$, $Y'_{t_k}=(0,\beta_0)$.
Ф.р. времён $\xi$ и $\xi'$ до прихода следующей заявки для процессов $Y_t$ и $Y_t'$ соответственно имеют вид:
\begin{eqnarray*}
\mathsf{P} \{\xi\leqslant s\}=F_1(s)=1-\exp\left(\displaystyle \int\limits_0^s\lambda(0,\alpha_0+u)\,\mathrm{d}\, u\right);
\\
\mathsf{P} \{\xi'\leqslant s\}=F_2(s)=1-\exp\left(\displaystyle \int\limits_0^s\lambda(0,\beta_0+u)\,\mathrm{d}\, u\right).
\end{eqnarray*}
В соответствии с Замечанием \ref{zv27} и (\ref{zv28}), общая часть этих распределений не меньше $\widehat\varkappa=\displaystyle  \frac{\lambda_0}{\Lambda}$.
Поэтому на пространстве $\displaystyle \prod\limits_{i=\ell_k+1} ^{\ell_k+3} (\Omega_i, \mathscr{F} _i , \mathsf{P} _i)$ (здесь $\nu_k=3$) в соответствии с Леммой \ref{zv29} строим $\xi$ и $\xi'$ таким образом, что $\mathsf{P} \{\xi=\xi'\}\geqslant \widehat\varkappa$.

В момент времени $t_{k+1}\stackrel{\text{\rm def}}{=\!=}  t_k+\min\{\xi,\xi'\}$ первая компонента хотя бы одного из процессов $Y_t$, $Y_t'$ возрастёт, и с вероятностью большей чем $\widehat\varkappa$ первые компоненты процессов $Y_t$, $Y_t'$ возрастут одновременно, т.е. произойдёт событие $ Y_{t_{k+1}}=Y_{t_{k+1}}'=(1,0;0)$.

\textbf{В. }Если случилось событие $Y_{t_{k}}=Y_{t_{k}}'=(1,0;0)$, при $t>t_k$ полагаем процессы $Y_t$ и $Y_t'$ равными, и строим процесс $Y_t$ (и равный ему $Y_t'$) по процедуре, описанной в пункте \textbf{А}.

Несложно заметить, что сконструированные по алгоритму \textbf{А-В} процессы $Y_t$ и $Y_t'$ удовлетворяют условиям (i)--(ii) Определения \ref{zv30}.

Покажем, что выполнено также условие (iii) Определения \ref{zv30}.

Рассмотрим процесс $Y_t'\stackrel{\mathscr{D}}{=}  X_t'$.
Если  $n_0'>0$, то пользуясь теми же соображениями, что и в Разделах \ref{zv20} и \ref{zv31}, можно показать, что остаточное время $\zeta_1$ первого периода занятости процесса $Y_t'$, зависящее от начального состояния $Y_0'$, имеет $k\in[0;K)$ конечных моментов.

Обозначим $\vartheta_0\stackrel{\text{\rm def}}{=\!=}  \zeta_0$, $\vartheta_k\stackrel{\text{\rm def}}{=\!=}  \inf\{t>\vartheta_k:\; Y'_t\in\mathscr{S} _0, Y'_{t-0}\notin\mathscr{S} _0\}$; $\vartheta_k$ -- это момент окончания $k$-го периода занятости процесса $Y_t'$.
В случае $Y_0'\in\mathscr{S} $ положим $\vartheta_0\stackrel{\text{\rm def}}{=\!=}  0$.

Поскольку распределение процесса $Y_t$ совпадает с распределением процесса $X_t$, верно неравенство   $\mathsf{P} \{X_{\vartheta_k}\in\mathscr{S} _0\}\geqslant e^{-\rho}$ --  см. (\ref{zv32}).

И, как мы уже видели, после совместного попадания процессов $Y_t$ и $Y_t'$ в множество $\mathscr{S} $ с вероятностью большей, чем $\widehat\varkappa$, при выходе из этого множества процессы ``склеиваются''.
Поэтому, снова обозначив $\tau\stackrel{\text{\rm def}}{=\!=} \inf\{t>0:\;Y_t=Y_t'\}$, а $\eta_k$ -- длина свободного периода $Y_t'$ после $\vartheta_k$, получаем, что
$ \displaystyle
\mathsf{P} \{\tau>\vartheta_k+\eta_k\}\leqslant (1-\varpi)^k.
$
Заметим, что $\tau=\tau(Y_0')=\tau(X_0')$ (поскольку $X_0=Y_0=(0,0)$ зафиксировано).

Очевидно, $\vartheta_k+\eta_k=\vartheta_0+\eta_0+\displaystyle  \sum\limits _{i=1}^k ( \zeta_i+\eta_i)$, где $\zeta_i$ -- $i$-й период занятости, а $\eta_i$ -- следующий за ним свободный период.

Обозначим $\mathscr{E} _k\stackrel{\text{\rm def}}{=\!=} \{Y_{\vartheta_k}\neq Y_{\vartheta_k}'\;\&\;Y_{\vartheta_k+\eta_k}= Y_{\vartheta_k+\eta_k}' \}$ -- событие, заключающееся в том, что склеивание произошло сразу после $k$-го попадания $Y_t'$ в множество $\mathscr{S} _0$.

\begin{proposition} \label{zv34} 
При $r\in[1; K)$ выполнено 
$\displaystyle
\mathbb{E} \,\tau^r\leqslant C_1(\varpi,r) \mathbb{E} \, \vartheta_0^r+ C_1(\varpi,r)\mathbb{E} \,\zeta_1^r+ C_3 (\varpi,r,\lambda_0,\Lambda),
$
константы $C_1$, $C_2$, $C_3$ вычисляются по приведённым ниже формулам.
\hfill \ensuremath{\triangleright}
\end{proposition}

\textbf{Доказательство.}
Сначала оценим $\mathbb{E} \,(\eta_k)^r|_{\mathscr{E} _k}$.
В случае события $\mathscr{E} _k$ сл.в. $\eta_k$ (в соответствии с Леммой \ref{zv29}) имеет плотность распределения
\begin{multline*}
  \varphi(s)= \min\left\{\lambda (Y_{\vartheta_k}+s)\exp \left(-\int\limits_0^ s\lambda(Y_{ \vartheta_k}+u)\,\mathrm{d}\,  u\right),\right.
  \\
\left.\lambda (Y'_{ \vartheta_k}+s)\exp \left(-\int\limits_0^s \lambda(Y'_{ \vartheta_k}+u)\right) \right\}
   \leqslant\Lambda e^{-\lambda_0 s}.
\end{multline*}
Поэтому $\mathbb{E} \,(\eta_k)^r|_{\mathscr{E} _k}\leqslant\displaystyle \frac{r! \Lambda}{\lambda_0^{r+1}}$.
Легко видеть, что и $\mathbb{E} \,(\eta_k)^r|_{\overline{\mathscr{E} }_k}\leqslant\displaystyle \frac{r! \Lambda}{\lambda_0^{r+1}}$; $\zeta_k$ от $\mathscr{E} _j$ не зависят при $k\neq j$.

Оценим $\mathbb{E} \,\tau^r$ (используя неравенство Йенсена в форме (\ref{zv33})):
\begin{multline*}
  \mathbb{E} \,\tau^r\leqslant \displaystyle  \sum\limits _{k=1}^\infty \mathbb{E} \left(\vartheta_0+\eta_0+\displaystyle  \sum\limits _{i=1}^k ( \zeta_i+\eta_i)\right)^r\mathbf{1} (\mathscr{E} _k)\leqslant
  \\
  \leqslant\displaystyle  \sum\limits _{k=1}^\infty (1-\varpi)^{k-1}\varpi (2k+2)^{r-1}\left( \left(\mathbb{E} \,\vartheta_0^r+(k+1)\displaystyle \frac{r! \Lambda}{\lambda_0^{r+1}}+k\mathbb{E} \,\zeta_1^r \right)\right)<
  \\
  <2^{r-1}\left(\mathbb{E} \,\vartheta_0^r\displaystyle  \sum\limits _{k=1}^\infty (1-\varpi)^{k-1}(k+1)^{r-1}+ \displaystyle \frac{r! \Lambda}{\lambda_0^{r+1}} \displaystyle  \sum\limits _{k=1}^\infty (1-\varpi)^{k-1}(k+1)^r+\right.
  \\
  +\left.\mathbb{E} \,\zeta_1^r\displaystyle  \sum\limits _{k=1}^\infty (1-\varpi)^{k-1}(k+1)^{r-1}k\right).
\end{multline*}
Предложение \ref{zv34} доказано.
\hfill \ensuremath{\blacksquare}

Итак,сконструированная выше пара $\mathscr{Z} _t=(Y_t,Y_t')$ -- успешная склейка процессов $X_t$ и $X_t'$, и получена явная оценка времени склеивания:
\begin{equation}\label{zv39}  
  \mathbb{E} \,(\tau(X_0'))^r\leqslant \mathcal{C}(r,\lambda_0, \Lambda,\mathbb{E} \,\vartheta_0^r).
\end{equation}

То есть
\begin{equation}\label{zv40}  
|\mathsf{P} \{X_t\in S\}-\mathsf{P} \{X_t'\in S\}|\leqslant \displaystyle \frac{\mathcal{C}(r, \lambda_0,\Lambda,\mathbb{E} \,\vartheta_0^r)}{t^r},
\end{equation}
однако следующий шаг -- интегрирование (\ref{zv35}) мы произвести не можем, поскольку стационарное распределение СМО $M|G|\infty$ (процесс $\stackrel{\circ}{X}_t$) и, тем более, $G|G|\infty$ (процесс ${X}_t$), \emph{неизвестно}.

\section{Процесс $X_t$ как регенерирующий процесс}
Напомним определение регенерирующего процесса.
\begin{definition}\label{rege} 
Случайный процесс $\{W_t,\,t\geqslant0\}$, заданный на вероятностном пространстве $(\Omega, \mathscr F, \mathbf{P})$, измеримый относительно фильтрации $\{\mathscr{F}_t,\,t\geqslant 0\}$, и принимающий значения в измеримом пространстве состояний \linebreak$(\mathscr{W},\mathscr{B}(\mathscr{W}))$  называется \textit{регенерирующим процессом}, если существует такая последовательность Марковских по отношению к фильтрации $\mathscr{F}_t$ моментов $\{\theta_i\}_{i\in \mathbb N}$, таких, что:
\begin{enumerate}
\item  $W_{\theta_i}=W_{\theta_j}$ при всех $i,j\in\mathbb N$;
\item  Случайные элементы $\Xi_i \stackrel{{\rm def}}{=\!\!\!=} \{W_t, \, t\in[\theta_i,\theta_{i+1}]\}$ ($i\in\mathbb N$) одинаково распределены и независимы в совокупности.\hfill \ensuremath{\triangleright}
\end{enumerate}
\end{definition}
Назовём периодами регенерации промежутки $[\theta_i,\theta_{i+1}]$; очевидно, длины периодов регенерации $\chi_i\stackrel{\text{\rm def}}{=\!=}  \theta_{i+1}-\theta_i$ одинаково распределены и независимы в совокупности.
Обозначим $\Phi(s)\stackrel{\text{\rm def}}{=\!=}  \mathsf{P} \{\chi_i\leqslant s\}$ -- (неизвестная нам) ф.р. сл.в. $\chi_i$

Известно, что если $\mathbb{E} \,\chi_i< \infty$, то распределение регенерирующего процесса $W_t$ слабо сходится к стационарному (инвариантному) распределению.

Заметим, что для введённого в Разделе \ref{zv16} процесса \linebreak$X_t=\left(n_t, x^{(0)}_t; x^{(1)}_t,x^{(2)}_t,\ldots,x^{(n_t)}_t\right)$ периоды регенерации разделяются на период занятости и свободный период: $\chi_k=\zeta_k+\eta_k$.

Моменты периода занятости оцениваются с помощью (\ref{zv36}) и (\ref{zv25}).

Моменты свободных периодов $\eta_i$, следующих за периодами занятости $\zeta_i$, оценены в Предложении \ref{zv34}.
Кроме того, период регенерации заведомо больше своей части -- свободного периода, а плотность распределения свободного периода $\varphi(s)\geqslant \lambda_0 e^{-\Lambda s}$ (Замечание \ref{zv27}).

Таким образом, известны оценки
\begin{equation}\label{zv43} 
  \begin{array}{l}
     \mathbb{E} \,\chi_k^r\leqslant 2^{r-1}(\mathbb{E} \,\zeta_1^r+\mathbb{E} \,\eta_1^r)=\mathcal{M}_r,\qquad
     \mathbb{E} \,\chi_i>\displaystyle \frac{\lambda_0}{\Lambda^2}.
  \end{array}
\end{equation}

Марковские моменты $\{\theta_i\}_{i\in \mathbb N}$ образуют вложенный процесс восстановления \linebreak $N_t=\displaystyle  \sum\limits _{i=1}^\infty \mathbf{1} (\theta_i<t)$.

Рассмотрим \emph{перескок} процесса $N_t$, а именно Марковский процесс $B_t\stackrel{\text{\rm def}}{=\!=}  t-\theta_{N_t}$ (с пространством состояний $\mathbb{R}_+$), распределение которого $\mathfrak{P}_t$, как известно (\cite{GBS,GK}), слабо сходится к стационарному инвариантному распределению $\mathfrak{P}$ с ф.р.
\begin{equation}\label{zv37} 
  \widetilde{\Phi}(s)=\displaystyle \frac{\displaystyle \int\limits_0^s(1-\Phi(u))\,\mathrm{d}\, u}{\mathbb{E} \,\chi_i}.
\end{equation}
К этому же распределению (\ref{zv37}) слабо сходится распределение \emph{недоскока} $R_t\stackrel{\text{\rm def}}{=\!=}  \theta_{N_t+1}-t$ процесса восстановления $N_t$.

Поскольку процесс $X_t$ -- регенерирующий, его распределение полностью определяется временем, прошедшим с момента начала периода регенерации.
Т.е. $\mathscr{P} _t=\mathbb{F} (B_t)$.

Теперь рассмотрим процесс $X_t$ с начальным состоянием $(0,0)$ и процесс $X_t'$ с произвольным начальным состоянием $X_0$.
Этим процессам соответствуют вложенные процессы восстановления $N_t$ и $N_t'$, а процессам $N_t$ и $N_t'$ соответствуют процессы перескоков $B_t$ и $B_t'$ и процессы недоскоков $R_t$ и $R_t'$.

Сконструировав успешную склейку $\mathscr{Z} _t=(Y_t,Y'_t)$  для процессов $X_t$ и $X_t'$, одновременно конструируется успешная склейка $\mathscr{W} _t=(D_t,D_t')$ процессов $B_t$ и $B_t'$.
Обозначим $\mathfrak{P}_t$ и $\mathfrak{P}_t'$ -- распределение процессов $B_t$ и $B_t'$ соответственно ($\mathsf{P} \{B_t\in S\}=\mathfrak{P}_t(S)$ и  $\mathsf{P} \{B_t'\in S\}=\mathfrak{P}_t'(S)$).

Из (\ref{zv39}) и (\ref{zv40}) сразу следует

\begin{equation}\label{zv41} 
  |\mathsf{P} \{B_t\in S\}-\mathsf{P} \{B_t'\in S\}|\leqslant \displaystyle \frac{\mathcal{C}(\varpi,r, \lambda_0,\Lambda,\mathbb{E} \,\vartheta_0^r)}{t^r},
\end{equation}
где $\vartheta_0=R_0'$ -- недоскок процесса $N_t'$ в момент времени $t=0$.
\begin{proposition}\label{zv42} 
\begin{multline*}
\|\mathfrak{P}_t-\mathfrak{P}\|_{TV}\leqslant
\\
\leqslant \frac{2\left( C_1(\varpi,r) \displaystyle \frac{\mathbb{E} \, \chi_i^{r+1}}{(r+1)\mathbb{E} \,\chi_i}+ C_1(\varpi,r)\mathbb{E} \,\zeta_1^r+ C_3 (\varpi,r,\lambda_0,\Lambda)\right)} {t^r} = \displaystyle \frac{\mathcal{C}_1(r,\lambda_0,\Lambda)}{t^r},
\end{multline*}
где $r\in[1,K-1)$.
\hfill \ensuremath{\triangleright}
\end{proposition}
\textbf{Доказательство.}
  Рассмотрим процесс перескока $B_t$ с начальным условием $B_0=0$ и процесс перескока $\widetilde{B}_t$, имеющий стационарное начальное распределение $\mathfrak{P}$; при всех $t\geqslant 0$
\begin{equation}\label{zv44} 
  \mathsf{P} \{\widetilde{B}_t\in S\}=\mathfrak{P}(S)=\int\limits_{S}\frac{1-\Phi(u)\,\mathrm{d}\, u}{\mathbb{E} \,\chi_i}.
\end{equation}

Подставляя в (\ref{zv41}) вместо $\vartheta_0=R_0'$ стационарный недоскок $\widetilde{R}_0$ с распределением (\ref{zv37}), получаем утверждение Предложения \ref{zv42}.

Учитывая (\ref{zv43}), можно вычислить верхнюю оценку величины $\mathcal{C}_1(r,\lambda_0,\Lambda)$.
\hfill \ensuremath{\blacksquare}

\begin{theorem} В условиях (\ref{zv10}) верно неравенство 
$\displaystyle
\|\mathscr{P} _t-\mathscr{P} \|_{TV}\leqslant \displaystyle \frac{\mathcal{C}_1(r,\lambda_0,\Lambda)}{t^r},
$
где $r\in[1,K-1)$.
\hfill \ensuremath{\triangleright}
\end{theorem}
\textbf{Доказательство.}
Как уже говорилось,  $\mathscr{P} _t=\mathbb{F}(\mathfrak{P}_t)$.
Соответственно, $\mathscr{P} =\mathbb{F}(\mathfrak{P})$.
Напомним, что $\chi_i$ -- длина периода регенерации, включающего в себя абсолютно непрерывные величины, т.е. распределение сл.в. $\chi_i$ имеет плотность.
Поэтому имеет плотность $\psi_t(s)$ перескок $B_t$; как уже говорилось (\ref{zv44}), известна ф.р. стационарного перескока, и её плотность равна $\psi(s)=\displaystyle \frac{1-\Phi(u)}{\mathbb{E} \,\chi_i}$.

Обозначим  $\sigma_t=\{s: \psi_t(s)> \psi(s)\}$.

Оценим
\begin{multline*}
  |\mathscr{P} _t(S)-\mathscr{P} (S)|=
  \\
  =\left|\int\limits_0^\infty \mathsf{P} \{X_t\in S|B_t=u\}\psi_t(u)\,\mathrm{d}\, u - \int\limits_0^\infty \mathsf{P} \{X_t\in S|B_t=u\}\psi(u)\,\mathrm{d}\, u\right|\leqslant
  \\
  \leqslant \max \left\{\int\limits_{\sigma_t}(\psi_t(u)-\psi (u))\,\mathrm{d}\, u, \int\limits_{\mathbb{R}_+\setminus \sigma_t}(\psi(u)-\psi_t (u))\,\mathrm{d}\, u \right\}=
  \\
  =\max\left\{|\mathfrak{P}_t(\sigma_t)- \mathfrak{P}(\sigma_t))|, |\mathfrak{P}_t(\mathbb{R}_+\setminus\sigma_t)- \mathfrak{P}(\mathbb{R}_+ \setminus\sigma_t)|\right\} \leqslant
  \\
  \leqslant \sup\limits _{\sigma\in \mathscr{B} (\mathbb{R})}| \mathfrak{P}_t(\sigma) -\mathfrak{P}(\sigma)|,
\end{multline*}
откуда $\|\mathscr{P} _t-\mathscr{P} \|_{TV}\leqslant \|\mathfrak{P}_t-\mathfrak{P}\|_{TV}$, что завершает доказательство.
\hfill \ensuremath{\blacksquare}

\begin{remark}\label{fin} 
Полученная оценка неоптимальна по следующим причинам.
\begin{enumerate}
  \item Применение метода склеивания для процессов в непрерывном времени, как правило, даёт неоптимальную оценку.
  \item При конструировании успешной склейки использовались только моменты попадания процесса $X_t'$ в свободное состояние; также можно было бы учесть и моменты попадания в свободное состояние процесса $X_t$.
  \item Условия (\ref{zv10}) не позволяют аккуратно оценить получаемые в ходе вычислений величины ($\mathcal{C}$ и $\mathcal{C}_1$).
      При исследовании конкретных прикладных задач эти величины могут быть оценены точнее.
\end{enumerate}
Тем не менее, полученная оценка может быть полезна в практических приложениях, поскольку имитационное моделирование не всегда способно дать достаточно качественные оценки.
\hfill \ensuremath{\triangleright}
\end{remark}

\paragraph{Благодарности.} Автор благодарит Л.Г.Афанасьеву, А.Ю.Веретенникова и А.Д.Маниту  за важные и полезные замечания при обсуждении представленной работы. Также автор благодарит Э.Ю.Калимулину за большую помощь в подготовке текста.
Работа выполнена при финансовой поддержке РФФИ (проект № 17-01-00633 А).